\documentclass[reqno,11pt]{amsart}
\usepackage{amsmath, amssymb, amsthm,amsfonts,pgfplotstable} 
\usepackage[english]{babel}
\usepackage{bbm}
\usepackage{graphicx}
\usepackage{url}
\usepackage{subcaption,caption}
\usepackage{pgfplots} \pgfplotsset{compat=newest}
\usepgfplotslibrary{groupplots}
\usepackage{epstopdf}
\usepackage[a4paper,bindingoffset=0.5cm,left=2.5cm,right=2.5cm,top=2.5cm,bottom=2cm,footskip=.8cm]{geometry}

\usepackage{rotating}
\usepackage{amsbsy,enumerate}

\usepackage{comment}
\usepackage{mathrsfs} 

\usepgfplotslibrary{groupplots}
\usepackage{tikz}
\usetikzlibrary{arrows.meta}
\usepackage{rotating}
\usepackage{amsbsy,enumerate}
\usepackage{graphicx}
\usepackage{comment}
\usepackage{mathrsfs}

\DeclareMathOperator{\E}{\mathbb{E}}

\newcommand{\s}{^\star}

\newcommand{\bs}{\boldsymbol}

\newcommand{\vb}{\vspace{3.2mm}}

\newcommand{\vt}{\vartheta}

\newcommand{\D}{\Delta}

\newcommand{\Dt}{\Delta t}
\renewcommand{\L}{\Lambda}
\newcommand{\p}{\partial}

\newcommand{\bvt}{\bm{\vartheta}}

\usepackage{color}

\newcommand{\ROOD}[1]{{\textcolor{black}{#1}}}
\newcommand{\RED}[1]{{\textcolor{black}{#1}}}

\newtheorem{remark}{Remark}
\newtheorem{example}{Example}

\newtheorem{proposition}{Proposition}

\allowdisplaybreaks

\begin{document}

\title[Linear \ROOD{stochastic fluid} networks]{Linear \ROOD{stochastic fluid} networks: \\rare-event simulation and Markov modulation}

\author{O.J. Boxma, E.J. Cahen, D. Koops, M. Mandjes}

\begin{abstract}
We consider a linear stochastic fluid network under Markov modulation, with a focus on the probability that the joint storage level attains a value in a rare set at a given point in time. The main objective is to develop efficient importance sampling algorithms with provable performance guarantees. For linear stochastic fluid networks without modulation, we prove that the number of runs needed (so as to obtain an estimate with a given precision) increases polynomially (whereas the probability under consideration decays essentially exponentially); for networks \RED{operating in the slow modulation regime,} our algorithm is asymptotically efficient. \ROOD{Our techniques are in the tradition of the rare-event simulation procedures that were developed for the sample-mean of i.i.d.\ one-dimensional light-tailed random variables, and intensively use the idea of exponential twisting. In passing, we also point out how to set up a recursion to evaluate the (transient and stationary) moments of the joint storage level in Markov-modulated linear stochastic fluid networks.}

\vspace{3mm}

\noindent {\sc Affiliations.}
O.\ Boxma is with E{\sc urandom} and Department of Mathematics and Computer Science, Eindhoven University of Technology, PO Box 513, 
5600 MB  Eindhoven, the Netherlands. 

\noindent
E.\ Cahen is with CWI, Science Park 123, 1098 XG Amsterdam, the Netherlands. 

\noindent D. Koops and M. Mandjes are with Korteweg-de Vries Institute for Mathematics, University of Amsterdam, Science Park 105-107, 1098 XG Amsterdam, the Netherlands. 

\noindent M. Mandjes is corresponding author. Email: {\scriptsize\tt m.r.h.mandjes@uva.nl}. Telephone: {\tiny  +} 31 20 525 5164. Fax: {\tiny  +} 31 20 525 7820.

\vspace{3mm}

\noindent {\sc Acknowledgments.} The research of O.\ Boxma, D.\ Koops and M.\ Mandjes was partly funded by the NWO Gravitation Project N{\sc etworks}, Grant Number 024.002.003. The research of O.\ Boxma was also partly funded by the Belgian Government, via the IAP Bestcom project. The research of E. Cahen was funded by an NWO grant, Grant Number 613.001.352.
\end{abstract}

\maketitle

\newpage
\section{Introduction}
\ROOD{\it Linear stochastic fluid networks}, as introduced in \cite{KW}, can be informally described as follows. Consider a network consisting of $L$ stations. Jobs, whose sizes are i.i.d.\ samples from some general $L$-dimensional distribution, arrive at the stations according to a Poisson process. At each of the nodes, in between arrivals the storage level decreases exponentially. Processed traffic is either transferred to the other nodes or leaves the network (according to a given routing matrix). In addition to this basic version of the linear stochastic fluid network, there is also its {\it Markov modulated} counterpart \cite{KS}, in which the arrival rate, the distribution of the job sizes, and the routing matrix depend on the state of an external, \ROOD{autonomously} evolving finite-state continuous-time Markov chain (usually referred to as the {\it background process}). 

Linear stochastic fluid networks can be seen as natural fluid counterparts of corresponding infinite-server queues. As such, they inherit several nice properties of those infinite-server queues. In particular, separate infinitesimally small fluid particles, moving through the network, do not interfere, and are therefore mutually independent. Essentially due to this property, linear stochastic fluid networks allow explicit analysis; in particular, the joint Laplace transform of the storage levels at a given point in time can be expressed in closed form as a function of the arrival rate, the Laplace transform of the job sizes and the routing matrix \cite[Thm. 5.1]{KW}. 

When Markov modulation is imposed, the analysis becomes substantially harder. Conditional on the path of the background process, again explicit expressions can be derived, cf.\ \cite[Thm.~1]{KS}. Unconditioning, however, cannot be done in a straightforward manner. As a consequence the results found are substantially less explicit than for the non-modulated linear stochastic fluid  network. In \cite{KS} also a system of ordinary differential equations has been set up that provides the transform of the stationary storage level; in addition, conditions are identified that guarantee the existence of such a stationary distribution.

\vb

In this paper we focus on rare events for Markov-modulated linear stochastic fluid networks. More specifically, in a particular scaling regime (parameterized by $n$) we analyze the probability $p_n$ that at a given point in time  the network storage vector is in a given rare set. By scaling the arrival rate as well as the rare set (which amounts to multiplying them by a scaling parameter $n$), the event of interest becomes increasingly rare. More specifically, under a Cram\'er-type assumption on the job-size distribution, application of large-deviations theory yields that $p_n$ decays (roughly) exponentially. As $p_n$ can be characterized only asymptotically, one could consider the option of using simulation to obtain precise estimates. The effectiveness, however, of such an approach is limited due to the rarity of the event under consideration: in order to get a reliable estimate, one needs sufficiently many runs in which the event occurs. This is the reason why one often resorts to simulation using {\it importance sampling} (or: {\it change of measure}). This is a variance reduction technique in \ROOD{which one replaces the actual probability measure by an alternative measure  under which the event under consideration is {\it not} rare}; correcting the simulation output with appropriate likelihood ratios yields an unbiased estimate. 

The crucial issue when setting up an importance sampling procedure concerns the choice of the alternative measure: one would like to select one that provides a substantial variance reduction, or is even (in some sense) optimal. The objective of this paper is to develop a change of measure which performs provably optimally.

Our  ultimate goal is to obtain an efficient simulation procedure for Markov-modulated linear stochastic fluid networks. We do so by (i)~first considering a single node without modulation, (ii) then  multi-node systems, still without modulation, and (iii) finally modulated multi-node systems. There are two reasons for this step-by-step setup:

\begin{itemize}
\item[$\circ$] For the non-modulated models we have more refined results than for the modulated models. More specifically, for the non-modulated models we have developed estimates for the number of runs $\Sigma_n$ required to obtain an estimate with predefined precision (showing that $\Sigma_n$ grows polynomially in the rarity parameter $n$), whereas for modulated models we can just prove that $\Sigma_n$ grows subexponentially. 
\item[$\circ$] In addition, this approach allows the reader to get gradually familiar with the concepts used in this paper. 
\end{itemize}
\RED{The construction and analysis of our importance sampling methodology is based on the ideas developed in  \cite{BM}; there the focus was on
addressing similar issues for a single-node Markov modulated infinite-server system. In line with \cite{BM}, we consider the regime in which the background process is `slow': while we (linearly) speed up the driving Poisson process, we leave the rates of the Markovian background process unalterned. }

\ROOD{A traditional, thoroughly examined, importance sampling problem concerns the sample mean $S_n$ of $n$ i.i.d.\ light-tailed random variables $X_1,\ldots,X_n$; the objective there is to estimate ${\mathbb P}(S_n \geqslant a)$ for $a> {\mathbb E} X_1$ and $n$ large. As described in \cite[Section VI.2]{AG}, in this situation importance sampling (i.e., sampling under an alternative measure, and translating the simulation output back by applying appropriate likelihood ratios) works extremely well. To this end, the distribution of the $X_i$\,s should be {\it exponentially twisted}. As it turns out, in our setup, the probability of our interest can be cast in terms of  this problem. Compared to the standard setup of sample means of one-dimensional  random variables, however, there are a few complications: (i) in our case it is not  a priori clear how to sample from the exponentially twisted distributions, (ii) we consider multi-dimensional distributions (i.e., rare-event probabilities that concern the storage levels of all individual buffers in the network), (iii) we impose Markov modulation. We refer to e.g.\ \cite{GJ,KMT} for earlier work on similar problems.}

In passing, we also point out how to set up a recursion to evaluate the (transient and stationary) moments of the joint storage level in Markov-modulated linear stochastic fluid networks (where the results in \cite{KS} are restricted to just the first two stationary moments at epochs that the background process jumps). 

\vb

The single-node model without modulation falls in the class of (one-dimensional) {\it shot-noise} models, for which efficient rare-event simulation techniques have been developed over the past, say, two decades. Asmussen and Nielsen \cite{AN} and Ganesh {\it et al.}\  \cite{GAN} consider the probability that a shot-noise process decreased by a linear drift ever exceeds some given level. 
Relying on sample-path large deviations results, an asymptotically efficient importance sampling algorithm is developed, under the same scaling as the one we consider in our paper. The major difference with our model (apart from the fact that we deal with considerably more general models, as we focus on networks and allow modulation) is that we focus on a rare-event probability that relates to the position of the process at a fixed point in time; in this setting we succeed in finding accurate estimates of the number of runs needed to get an estimate of given precision. 

\ROOD{There is a vast body of literature related to the broader area of rare-event simulation for queueing systems. We refer to the literature overviews \cite{BLM,JS}; interesting recent papers include \cite{AK, CMZ, SEZ}. }

\vb

This paper is organized as follows. In Section \ref{S2} the focus is on a single-node network, without Markov modulation (addressing complication (i) above), Section \ref{S3} addresses the extension to multi-node systems (addressing complication (ii)), and in Section~\ref{S4} the feature of modulation is added (addressing complication (iii)). In each of these three sections, we propose a change of measure, quantify its performance, and demonstrate its efficiency through simulation experiments. \RED{In Section \ref{S41} we include the explicit expressions for the moments  in Markov-modulated linear stochastic fluid networks.}
A discussion and concluding remarks are found in Section \ref{S5}.

\section{Single resource, no modulation}\label{S2}
To introduce the concepts we work with in this paper, we analyze in this section  a linear stochastic fluid network consisting of a single node, in which  the input is just compound Poisson (so no Markov modulation is imposed). More precisely, in the model considered, 
jobs arrive according to a Poisson process with rate $\lambda$, bring along i.i.d.\ amounts of work (represented by the sequence of i.i.d.\ random variables $(B_1,B_2,\ldots)$), and the workload level decays exponentially at a rate $r>0.$ This model belongs to the class of {\it shot-noise processes}. As mentioned in the introduction, we gradually extend the model in the next sections. 
\subsection{Preliminaries}
We first present a compact representation for the amount of work in the system at time $t$, which we denote by $X(t)$, through its moment generating function. To this end, let $N(t)$ denote a Poisson random variable with mean $\lambda t$, and $(U_1,U_2,\ldots)$ i.i.d.\ uniformly distributed random variables (on the interval $[0,t]$). Assume in addition that the random objects $(B_1,B_2,\ldots)$, $N(t)$, and $(U_1,U_2,\ldots)$ are independent. Then it is well-known that the value of our shot-noise process at time $t$ can be expressed as
\begin{equation}\label{swap}X(t) = \sum_{j=1}^{N(t)} B_j {\rm e}^{-r(t-U_j)}\stackrel{\rm d}{=}\sum_{j=1}^{N(t)} B_j {\rm e}^{-rU_j},\end{equation}
where the distributional equality is a consequence of the fact that the distribution of $U$ is symmetric  on the interval $[0,t]$. 
It is easy to compute the moment generating function ({\sc mgf}) of $X(t)$, by conditioning on the value of $N(t)$:
\begin{eqnarray}\nonumber M(\vt):={\mathbb E}\,{\rm e}^{\vartheta X(t)}&=&
\sum_{k=0}^\infty
{\rm e}^{-\lambda t}\frac{(\lambda t)^k}{k!} \left({\mathbb E}\exp (\vt B \,{\rm e}^{-rU})\right)^k
\\&=&\label{mgf}
\exp\left(\lambda\int_0^t \left(\beta( {\rm e}^{-ru}\,\vt)-1\right){\rm d}u\right),\end{eqnarray}
where $\beta(\cdot)$ is the {\sc mgf}\:corresponding to $B$ (throughout assumed to exist).
By differentiating and inserting $\vt=0$, it follows immediately that
\[{\mathbb E}\,X(t) = \frac{\lambda}{r}(1-{\rm e}^{-rt})\, {\mathbb E}\,B=: m(t).\]
Higher moments can be found by repeated differentiation. We note that, as $t$ is held fixed throughout the document, we often write $N$ rather than $N(t).$

\subsection{Tail probabilities, change of measure}
The next objective is to consider the asymptotics of the random variable $X(t)$ under a particular scaling. In this scaling we let the arrival rate be $n\lambda$ rather than just $\lambda$, for ${n\in{\mathbb N}}$. The value of the shot-noise process is now given by
\[Y_n(t) := \sum_{i=1}^n X_i(t),\]
with the vector  $(X_1(t),\ldots,X_n(t))$ consisting of i.i.d.\ copies of the random variable $X(t)$ introduced above; here the infinite divisibility of a Compound Poisson distribution is used. 

Our goal is to devise techniques to analyze the tail distribution of $Y_n(t).$
Standard theory now provides us with the asymptotics of 
\[p_n(a)={\mathbb P}(Y_n(t) \geqslant na)\] for some $a>m(t)$; we are in the classical `Cram\'er setting' \cite[Section 2.2]{DZ} if it is assumed that $M(\vt)$ is finite in a neighborhood around the origin (which requires that the same property is satisfied by $\beta(\cdot)$).  \ROOD{Let $I(a)$ and $\vartheta\s\equiv \vt\s(a)$, respectively, be defined as 
\[I(a):=\sup_\vartheta\big(\vartheta a - \log M(\vt)\big),\:\:\:\:
\vartheta\s:= \arg\sup_\vartheta\big(\vartheta a - \log M(\vt)\big),\]
with $M(\cdot)$ as above.
Using `Cram\'er', we obtain that, under mild conditions, 
\[\lim_{n\to\infty}\frac{1}{n}\log p_n(a) =-I(a)= -\vartheta\s a +\log M(\vt\s).\]
More refined asymptotics are available as well; we get back to this issue in Section \ref{EFF}. }

As these results apply in the regime that $n$ is large, a relevant issue concerns the development of efficient techniques to estimate $p_n(a)$ through simulation. An important rare-event simulation technique is  importance sampling, relying on the commonly used exponential twisting technique. We now investigate how to construct the exponentially twisted version ${\mathbb Q}$ (with twist $\vt\s$) of the original probability measure ${\mathbb P}.$ The main idea is that under ${\mathbb Q}$ the $X_i(t)$ have mean $a$, such that under the new measure the event under study is not rare anymore.

More concretely, exponential twisting with parameter $\vt\s$ means that under the new measure ${\mathbb Q}$, the $X_i(t)$ should have the {\sc mgf}
\begin{equation}\label{newMGF} {\mathbb E}_{\mathbb Q}\, {\rm e}^{\vt X(t)} = \frac{{\mathbb E}\, {\rm e}^{(\vt+\vt\s) X(t)} }{{\mathbb E}\, {\rm e}^{\vt\s X(t)} }
=\frac{M(\vt+\vt\s)}{M(\vt\s)};\end{equation}
under this choice the random variable has the desired mean:
\[{\mathbb E}_{\mathbb Q}\, {X(t)} = \frac{M'(\vt\s)}{M(\vt\s)} = a.\]
The question is now: how to sample a random variable that has this {\sc mgf}? To this end, notice that $M(\vt) =\exp(-\lambda t+\lambda t\,{\mathbb E} \exp(\vt B {\rm e}^{-rU}))$ and
\[ M(\vt+\vt\s)   =\sum_{k=0}^\infty 
{\rm e}^{-\lambda t}\frac{(\lambda t\, {\mathbb E}\exp(\vt\s B\,{\rm e}^{-rU}))^k}{k!} \left(\frac{{\mathbb E}\exp((\vt+\vt\s) B\,{\rm e}^{-rU})}{{\mathbb E}\exp(\vt\s B\,{\rm e}^{-rU})}\right)^k,\]
such that (\ref{newMGF}) equals 
\[\sum_{k=0}^\infty 
\exp(-\lambda t\, {\mathbb E}\exp(\vt\s B\,{\rm e}^{-rU}))\frac{(\lambda t\, {\mathbb E}\exp(\vt\s B\,{\rm e}^{-rU}))^k}{k!} \left(\frac{{\mathbb E}\exp((\vt+\vt\s) B\,{\rm e}^{-rU})}{{\mathbb E}\exp(\vt\s B\,{\rm e}^{-rU})}\right)^k .\]
From this expression we can see how to sample the $X_i(t)$ under ${\mathbb Q}$, as follows. 
In the first place we conclude that under ${\mathbb Q}$ the number of arrivals becomes Poisson with mean
\begin{equation}
\label{PPP}\lambda t\, {\mathbb E}\exp(\vt\s B\,{\rm e}^{-rU}) =\lambda \int_0^t  \beta({\rm e}^{-ru}\,\vt\s){\rm d}u,\end{equation}
rather than $\lambda t$ (which is an increase). Likewise, it entails that under ${\mathbb Q}$ the distribution of the $B_j {\rm e}^{-rU_j}$ should be twisted by $\vt\s$, in the sense that these random variables should have under ${\mathbb Q}$ the {\sc mgf}\:
\[{\mathbb E}_{\mathbb Q}\exp((\vt+\vt\s) B\,{\rm e}^{-rU})=\frac{{\mathbb E}\exp((\vt+\vt\s) B\,{\rm e}^{-rU})}{{\mathbb E}\exp(\vt\s B\,{\rm e}^{-rU})}.\] We now point out how such a random variable should be sampled. To this end, observe that
\[{\mathbb E}\exp((\vt+\vt\s) B\,{\rm e}^{-rU}) = \int_0^t
\frac{\beta({\rm e}^{-ru}(\vt+\vt\s))}{\beta( {\rm e}^{-ru}\,\vt\s)} \frac{1}{t}\,\beta( {\rm e}^{-ru}\,\vt\s){\rm d}u,\]
so that
\[{\mathbb E}_{\mathbb Q}\exp((\vt+\vt\s) B\,{\rm e}^{-rU})
=\int_0^t
\frac{\beta({\rm e}^{-ru}(\vt+\vt\s))}{\beta( {\rm e}^{-ru}\,\vt\s)} \frac{\beta( {\rm e}^{-ru}\,\vt\s)}{\displaystyle \int_0^t \beta( {\rm e}^{-rv}\,\vt\s) {\rm d}v}{\rm d}u.\]
From this representation two conclusions can be drawn. In the first place, supposing there are $k$ arrivals, then the arrival epochs $U_1,\ldots,U_k$
are i.i.d.\ under ${\mathbb Q}$, with the density given by
\[f_U^{\mathbb Q}(u) = \frac{\beta( {\rm e}^{-ru}\,\vt\s)}{\displaystyle \int_0^t \beta({\rm e}^{-rv}\,\vt\s )\,{\rm d}v}.\]
In the second place, given that the $k$-th arrival occurs at time $u$,  the density of the corresponding job size $B_k$ should be exponentially twisted by ${\rm e}^{-ru}\,\vt\s$ (where each of the job sizes is sampled independently of everything else). 

Now that we know how to sample from ${\mathbb Q}$ it is straightforward to implement the importance sampling. Before we describe its complexity (in terms of the number of runs required to obtain an estimate with given precision), we first provide an example in which we demonstrate how the change of measure can be performed.

\begin{example}\label{Ex1}{\em In this example we consider the case that the $B_i$ are exponentially distributed with mean $\mu^{-1}$. Applying the transformation $w:=  {\rm e}^{-ru}\,\vt/\mu,$
it is first seen that
\begin{eqnarray*}
\int_0^s \beta( {\rm e}^{-ru}\,\vt) {\rm d}u &=&
\int_0^s \frac{\mu}{\mu- {\rm e}^{-ru}\,\vt\,}{\rm d}u=\frac{1}{r}\int_{{\rm e}^{-rs}\vt/\mu}^{\vt/\mu}
\frac{1}{1-w}\frac{1}{w}{\rm d}w\\
&=&\frac{1}{r}\left[\log\frac{w}{1-w}\right]_{ {\rm e}^{-rs}\vt/\mu}^{\vt/\mu}=\frac{1}{r}\log\left(\frac{\mu {\rm e}^{rs}-\vt}{\mu-\vt}\right).
\end{eqnarray*}
As $\vt\s$ solves the equation $M'(\vt\s)/M(\vt\s) = a$, we obtain the quadratic equation
\[m(t) = a\left(1-\frac{\vt}{\mu}\right)\left(1-\frac{\vt}{\mu}{\rm e}^{-rt}\right),\]
leading to
\[\vt\s = \frac{\mu {\rm e}^{rt}}{2}\left((1+{\rm e}^{-rt}) -\sqrt{(1-{\rm e}^{-rt})^2+4{\rm e}^{-rt}\frac{m(t)}{a}}\right)\]
(where it is readily checked that $\vt\s\in(0,\mu)$).

Now we compute what the alternative measure ${\mathbb Q}$ amounts to. In the first place, 
the number of arrivals should become Poisson with parameter
\[\frac{\lambda}{r}\log\left(\frac{\mu {\rm e}^{rt}-\vt\s}{\mu-\vt\s}\right)\]
(which is larger than $\lambda t$).
In addition, we can check that
\[F_U^{\mathbb Q}(u) :={\mathbb Q}(U\leqslant u) = \log\left(\frac{\mu {\rm e}^{ru}-\vt\s}{\mu-\vt\s}\right)\left/
\log\left(\frac{\mu {\rm e}^{rt}-\vt\s}{\mu-\vt\s}\right)\right.\]
(rather than $u/t$). 
The function $F_U^{\mathbb Q}(u)$ has the value $0$ for $u=0$ and the value $1$ for $u=t$, and is concave. This concavity reflects that the arrival epochs of the shots tend to be closer to $0$ under ${\mathbb Q}$ than under ${\mathbb P}$. \ROOD{This is because we identified each of the $U_i$ with $t$ minus the actual corresponding arrival epoch in (\ref{swap}); along the most likely path of $Y_n(t)$ itself the shots will be typically closer to $t$ under ${\mathbb Q}$.}
Observe that one can sample $U$ under ${\mathbb Q}$ using the classical inverse distribution function method \cite[Section II.2a]{AG}: with $H$ denoting a uniform number on $[0,1)$, we obtain such a sample by
\[\frac{1}{r}\log\left(\left({\rm e}^{rt}-\frac{\vt\s}{\mu}\right)^H\left(1-\frac{\vt\s}{\mu}\right)^{1-H}+\frac{\vt\s}{\mu}\right).\]
Also, conditional on a $U_i$ having attained the value $u$, the jobs $B_i$ should be sampled from an exponential distribution with mean $(\mu- {\rm e}^{-ru}\,\vt\s)^{-1}$. }\end{example}

\begin{remark}{\em 
In the model we study in this section, the input of the linear stochastic fluid network is a compound Poisson process. As pointed out in \cite{KW} the class of inputs can be extended to the more general class of increasing L\'evy processes in a straightforward manner. }
\end{remark}

\subsection{Efficiency properties of importance sampling procedure} 
\label{EFF} 
In this subsection we analyze the performance of the procedure introduced in the previous section. The focus is on a characterization of the number of runs needed to obtain an estimate with a given precision (at a given confidence level). 
 
In every run $Y_n(t)$ is sampled under ${\mathbb Q}$, as pointed out above. As ${\mathbb Q}$ is an implementation of an exponential twist (with twist $\vt\s$), the likelihood ratio (of sampling  $Y_n(t)$ under ${\mathbb P}$ relative to ${\mathbb Q}$) is given by
\[L= \frac{{\rm d}{\mathbb P}}{{\rm d}{\mathbb Q}} = {\rm e}^{-\vt\s Y_n(t)} {\rm e}^{n\,\log M(\vt\s)}.\]
\ROOD{In addition, define $I$ as the indicator function of the event $\{Y_n(t)\geqslant na\}$.} \ROOD{Clearly, ${\mathbb E}_{\mathbb Q}(LI)=p_n(a)$.} \ROOD{We keep generating samples $LI$ 
(under ${\mathbb Q}$), and estimate $p_n(a)$ by the corresponding sample mean, until the ratio of the half-width of the confidence interval (with critical value $T$) and the estimator drops below some predefined $\varepsilon$ (say, 10\%).}
Under ${\mathbb P}$ the number of runs needed is effectively inversely proportional to $p_n(a)$, hence exponentially increasing in $n$. We now focus on quantifying the reduction of the number of runs when using the  importance sampling procedure we described above, i.e., the one based on the measure ${\mathbb Q}$.

Using a Normal approximation, it is a standard reasoning that when performing $N$ runs the ratio of the half-width of the confidence interval and the estimator is approximately
\[\frac{1}{p_n(a)}\cdot \frac{T}{\sqrt{N}}\sqrt{{\mathbb V}{\rm ar}_{\mathbb Q} (L^2 I)},\]
and hence the number of runs needed is roughly
\[\Sigma_n :=\frac{T^2}{\varepsilon^2} \,\frac{{\mathbb V}{\rm ar}_{\mathbb Q} (L^2 I)}{(p_n(a))^2}.\]
We now analyze how $\Sigma_n$ behaves as a function of the `rarity parameter' $n$. Due to the Bahadur-Rao result \cite{BR}, with $f_n\sim g_n$ denoting $f_n/g_n\to 1$ as $n\to\infty$,
\begin{equation}
\label{FIRM}
p_n(a)={\mathbb E}_{\mathbb Q}(LI)\sim \frac{1}{\sqrt{n}} \frac{1}{\vt\s\sqrt{2\pi\tau}} {\rm e}^{-nI(a)},\:\:\:\tau:=\left.\frac{{\rm d}^2}{{\rm d}\vt^2}\log M(\vt)\right|_{\vt=\vt\s}.\end{equation}
Using the same proof technique as in \cite{BR}, it can be shown that 
\begin{equation}
\label{SECM}
{\mathbb E}_{\mathbb Q}(L^2I)  \sim \frac{1}{\sqrt{n}} \,\frac{1}{2\vt\s\sqrt{2\pi\tau}} {\rm e}^{-2nI(a)};\end{equation}
\RED{see Appendix A for the underlying computation. It also follows that ${\mathbb E}_{\mathbb Q}(L^2I)  \sim {\mathbb V}{\rm ar}_{\mathbb Q} (L^2 I)$.} 

We can use these asymptotics, to conclude that under ${\mathbb Q}$ the number of runs required grows slowly in $n$. More specifically,
$\Sigma_n$ is essentially proportional to $\sqrt{n}$ for $n$ large. This leads to the following result; \RED{cf.\ \cite[Section 2]{BGL} for related findings in a more general context.}
\begin{proposition}
As $n\to\infty$,
\begin{equation}\label{alpha0}\Sigma_n\sim \alpha\,\sqrt{n},\:\:\:\:\alpha=\frac{T^2}{\varepsilon^2} \, {\vt\s}\cdot\frac{1}{2}\,\sqrt{2\pi \tau}.\end{equation}
\end{proposition}

\subsection{Simulation experiments}\label{S24}
\begin{figure}
\centering
%
%
%
%
%
\includegraphics[width=.9\textwidth]{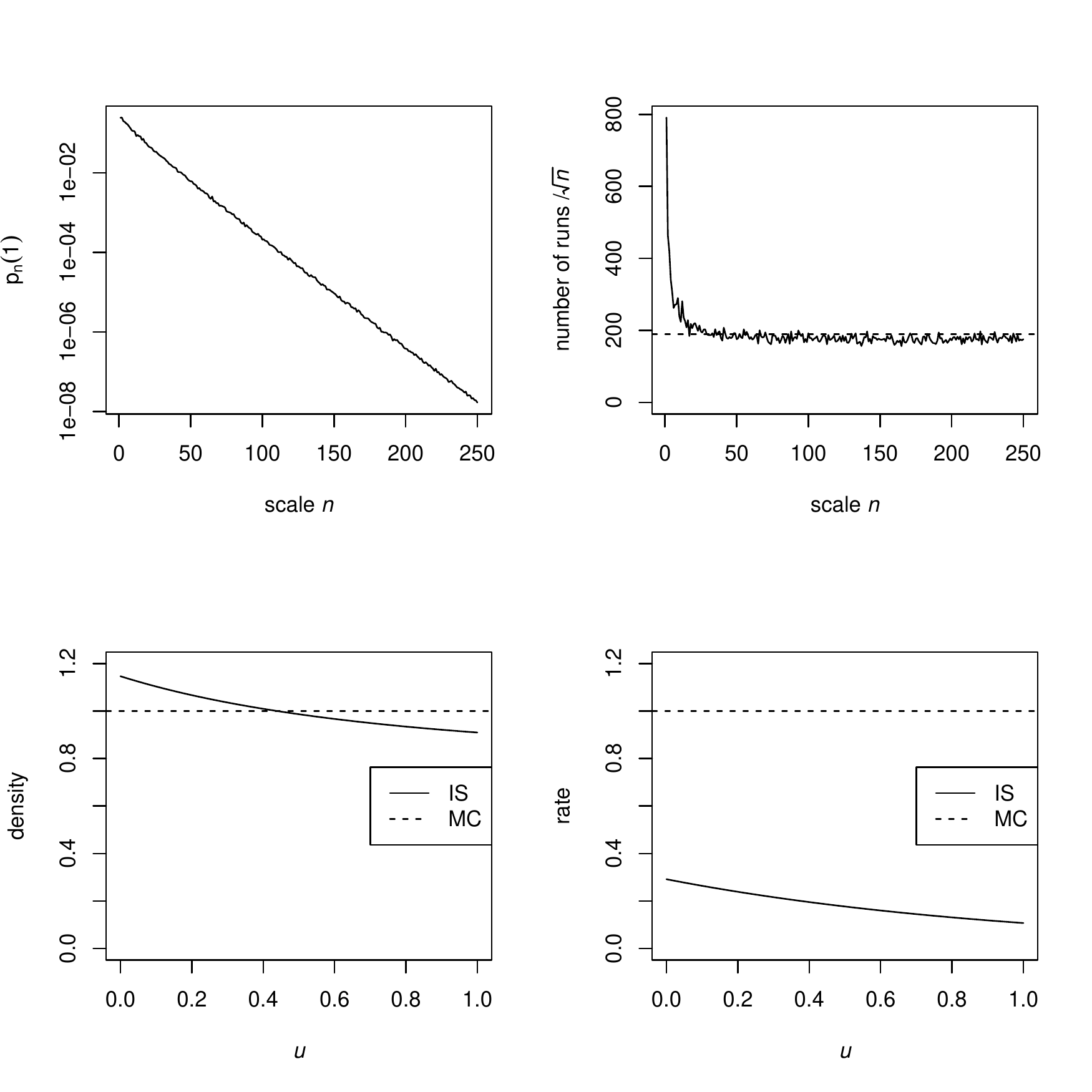}
\caption{Numerical results for Section \ref{S24}.}
\label{F1}
\end{figure}

In this subsection we present numerical results for the single-node model without Markov modulation. We focus on the case of exponential jobs, as in Example~\ref{Ex1}. We simulate until the estimate has reached the precision $\varepsilon = 0.1$, with confidence level $0.95$ (such that the critical value is $T=1.96$). The parameters chosen are: $t=1$, $r=1$, $\lambda = 1$, and $\mu=1$. We set $a=1$ (which is larger than $m(t) = 1-e^{-1}$). As it turns out, $\vt\s =0.2918$ and
\[\tau =\frac{\lambda}{r}\left(\frac{1}{(\mu-\vt\s)^2}-\frac{1}{(\mu{\rm e}^{rt}-\vt\s)^2}\right)= 1.8240.\]
The top-left panel of Fig.\ \ref{F1} confirms the exponential decay of the probability of interest, as a function of $n$. 
In the top-right panel we verify that the number of runs indeed grows proportionally to $\sqrt{n}$; 
the value of $\alpha$, as defined in (\ref{alpha0}), is $198.7$, which is depicted by the  horizontal line.
The bottom-left panel shows the density of the arrival epochs, which confirms that the arrival epochs tend to be closer to $0$ under ${\mathbb Q}$ than under ${\mathbb P}$; recall that under ${\mathbb P}$ these epochs are uniformly distributed on $[0,t]$. Recall that we reversed time in (\ref{swap}): for the actual shot-noise system that we are considering, it means that in order to reach the desired level at time $t$, the arrival epochs tend to be closer to $t$ under ${\mathbb Q}$ than under ${\mathbb P}$.
The bottom-right panel presents the rate of the exponential job sizes as a function of $u$. Using (\ref{PPP}), the arrival rate under ${\mathbb Q}$ turns out to be $1.2315$.

\section{Multi-node systems, no modulation}\label{S3}
In this section we consider multi-node stochastic fluid linear stochastic fluid networks, of the type analyzed in the work by Kella and Whitt \cite{KW}. It is instructive to first consider  the simplest multi-node system: a tandem network without external input in the downstream node \ROOD{and no traffic leaving after having been served by the upstream node} (and rate $r_\ell$ for node $\ell$, $\ell=1,2$); later we extend the ideas developed to general linear stochastic fluid networks. 

\subsection{Preliminaries}
As mentioned above, we first consider the two-node tandem.
The content of the first node is, as before,
\[X^{(1)}(t) =\sum_{j=1}^N B_j {\rm e}^{-r_1(t-U_j)}\]
(with $N$ having a Poisson distribution with mean $\lambda t$),
but it can be argued that the content of the second node satisfies a similar representation. More specifically, using the machinery developed in \cite{KW}, it turns out that
\begin{equation}
\label{distr}X^{(2)}(t) = \sum_{j=1}^N B_j \frac{r_1}{r_1-r_2} \left({\rm e}^{-r_2(t-U_j)}-{\rm e}^{-r_1(t-U_j)}\right)\stackrel{\rm d}{=}
 \sum_{j=1}^N B_j \frac{r_1}{r_1-r_2} \left({\rm e}^{-r_2 U_j}-{\rm e}^{-r_1 U_j}\right).\end{equation}
As before, perform the scaling by $n$, meaning that the arrival rate $\lambda$ is inflated by a factor $n$. It leads to the random vectors $(X^{(1)}_1(t),\ldots,X^{(1)}_n(t))$ and $(X^{(2)}_1(t),\ldots,X^{(2)}_n(t))$. With these vectors we can define $Y_n^{(1)}(t)$ and $Y_n^{(2)}(t)$, analogously to how this was done in the single-node case; these two random quantities represent the contents of the upstream resource and the downstream resource, respectively.

The state of this tandem system can be uniquely characterized in terms of its (bivariate) moment generating function. The technique to derive an explicit expression is by relying on the above distributional equality (\ref{distr}). Again, the key step is to condition on the number of shots that have arrived in the interval $[0,t]$: with ${\boldsymbol \vt} = (\vt_1,\vt_2)$,
\begin{eqnarray}\nonumber
M({\boldsymbol \vt})&:=& {\mathbb E}\, {\rm e}^{\vt_1 X^{(1)}(t)+\vt_2 X^{(2)}(t)}\\ \nonumber
&=&\sum_{k=0}^\infty {\rm e}^{-\lambda t}\frac{(\lambda t)^k}{k!}\left({\mathbb E}
\exp\left(\vt_1 B{\rm e}^{-r_1U} +\vt_2B\frac{r_1}{r_1-r_2}\left({\rm e}^{-r_2 U}-{\rm e}^{-r_1 U}\right)
\right)\right)^k\\ \nonumber
&=& \sum_{k=0}^\infty {\rm e}^{-\lambda t}\frac{(\lambda t)^k}{k!}\left(\int_0^t \frac{1}{t}{\mathbb E}
\exp\left(\vt_1 B{\rm e}^{-r_1u} +\vt_2B\frac{r_1}{r_1-r_2}\left({\rm e}^{-r_2 u}-{\rm e}^{-r_1 u}\right)
\right){\rm d}u\right)^k\\ \nonumber
&=& \sum_{k=0}^\infty {\rm e}^{-\lambda t}\frac{(\lambda t)^k}{k!}\left(\int_0^t \frac{1}{t}\beta\left( {\rm e}^{-r_1u}\vt_1 +\frac{r_1}{r_1-r_2}\left({\rm e}^{-r_2 u}-{\rm e}^{-r_1 u}\right)\vt_2
\right){\rm d}u\right)^k\\
&=&\exp\left( \lambda \int_0^t \left(\beta\left( {\rm e}^{-r_1u} \vt_1+\frac{r_1}{r_1-r_2}\left({\rm e}^{-r_2 u}-{\rm e}^{-r_1 u}\right)\vt_2
\right)-1\right){\rm d}u \right).\label{mgf2}
\end{eqnarray}

The above computation is for the two-node tandem system, but the underlying procedure can be extended to the case of networks with more than 2 nodes, and external input in each of the nodes. To this end, we consider the following  network consisting of $L$ nodes. Jobs are generated according to a Poisson process. \ROOD{At an arrival epoch, an amount  is added to the content of each of the resources $\ell\in\{1,\ldots,L\}$, where the amount added to resource $\ell$ is distributed as the (non-negative) random variable $B^{(\ell)}$;} $\beta({\boldsymbol \vt})$, with ${\bs\vt}\in{\mathbb R}^L$, is the joint {\sc mgf} of $B^{(1)}$ up to $B^{(L)}$ (note that the components are not assumed independent). In addition, let the traffic level at node $\ell$ decay exponentially with rate $r_\ell$
(i.e., the value of the output rate is linear in the current level, with proportionality constant $r_\ell$). A deterministic fraction $p_{\ell\ell'}\geqslant 0$ ($\ell\not=\ell'$) is then fed into node $\ell'$, whereas a fraction $p_{\ell\ell}\geqslant 0$ leaves the network (with $\sum_{\ell'=1}^Lp_{\ell\ell'} = 1$). We denote $r_{\ell\ell'} := r_\ell p_{\ell\ell'}.$ As an aside we mention that this general model covers models in which some arrivals (of the Poisson process with parameter $\lambda$) actually lead to arrivals at only a subset of the $L$ queues (since the job sizes $B^{(1)},\ldots,B^{(L)}$ are allowed to equal $0$).

We now point out how the joint buffer content process can be analyzed. Again our objective is to evaluate the moment generating function. 
Define the matrix $R$ as follows: its $(\ell,\ell)$-th entry is $r_{\ell\ell}+\sum_{\ell'\not=\ell}r_{\ell\ell'}$, whereas its $(\ell,\ell')$-th entry (with $\ell\not=\ell'$) is $-r_{\ell\ell'}$.
We have, according to Kella and Whitt \cite{KW}, with $N$ again Poisson with mean $\lambda t$, the following distributional equality: for any $\ell\in\{1,\ldots,L\}$,
\[X^{(\ell)}(t) = \sum_{\ell'=1}^L\sum_{j=1}^N B_j^{(\ell')} ({\rm e}^{-R(t-U_j)})_{\ell'\ell}.\]
It means we can compute the joint {\sc mgf} of $X^{(1)}(t)$ up to $X^{(L)}(t)$ as follows, cf.\ \cite[Thm. 5.1]{KW}:
\begin{eqnarray*}
M({\boldsymbol \vt})&:=& {\mathbb E}\, \exp\left(\sum_{\ell=1}^L\vt_\ell X^{(\ell)}(t)\right)\\
&=&\sum_{k=0}^\infty {\rm e}^{-\lambda t}\frac{(\lambda t)^k}{k!}\left({\mathbb E}
\exp\left(\sum_{\ell=1}^L\vt_\ell \sum_{\ell'=1}^L
B^{(\ell')} ({\rm e}^{-R(t-U)})_{\ell'\ell }
\right)\right)^k\\
&=& \sum_{k=0}^\infty {\rm e}^{-\lambda t}\frac{(\lambda t)^k}{k!}\left(\int_0^t \frac{1}{t}{\mathbb E}
\exp\left(\sum_{\ell=1}^L\vt_\ell \sum_{\ell'=1}^L B^{(\ell')} ({\rm e}^{-Ru})_{\ell'\ell}\right)
{\rm d}u\right)^k\\
&=& \sum_{k=0}^\infty {\rm e}^{-\lambda t}\frac{(\lambda t)^k}{k!}\left(\int_0^t \frac{1}{t}\beta\left(\sum_{\ell=1}^L
({\rm e}^{-Ru})_{1\ell}\vt_\ell,\ldots,\sum_{\ell=1}^L
({\rm e}^{-Ru})_{L\ell}\vt_\ell\right)
{\rm d}u\right)^k\\
&=&\exp\left(-\lambda t + \lambda \int_0^t\beta\left(\sum_{\ell=1}^L
({\rm e}^{-Ru})_{1\ell}\vt_\ell,\ldots,\sum_{\ell=1}^L
({\rm e}^{-Ru})_{L\ell}\vt_\ell\right){\rm d}u \right)\\
&=&\exp\left(\lambda \int_0^t\left(\beta\left({\rm e}^{-Ru}\,{\bs\vt}\right)-1\right){\rm d}u \right),
\end{eqnarray*}
which is the matrix/vector-counterpart of the expression (\ref{mgf}) that we found in the single-node case; for the two-node case the special form (\ref{mgf2}) applies.

\subsection{Tail probabilities, change of measure} \label{S32}
In this subsection we introduce the change of measure that we use in our importance sampling approach.
Many of the concepts are analogous to concepts used for the single-node case in Section 2.

Define (in self-evident notation)
\[p_n({\boldsymbol a}) :=
{\mathbb P}\left(Y^{(1)}_n(t)\geqslant na_1, \ldots,Y^{(L)}_n(t)\geqslant na_L\right).\]
Due to the multivariate version of Cram\'er's theorem, with $A:=[a_1,\infty)\times\cdots\times [a_L,\infty)$,
\begin{equation}\label{logass}\lim_{n\to\infty}\frac{1}{n}\log p_n({\boldsymbol a}) = -\inf_{{\bs b}\in A}I({\bs b}),\:\:\:\mbox{where}\:\:\:I({\bs b}):=\sup_{\bs\vt}\left(\langle{\bs \vt},{\bs b}\rangle  - \log M({\bs \vt})\right).\end{equation}
More refined asymptotics than the logarithmic asymptotics of (\ref{logass}) are available as well, but these are not yet relevant in the context of the present subsection; we return to these `exact asymptotics' in Section \ref{S33}.

We assume that the set $A$ is `rare', in the sense that 
\[{\bs m}(t)\not\in A,\:\:\:\mbox{with} \:\:\:m_i(t):=\left.\frac{\partial M({\bs \vt})}{\partial \vt_i}\right|_{{\bs \vt}={\bs 0}}.\]

Let us now construct the importance sampling measure. Let ${\boldsymbol \vartheta}\s$ be the optimizing ${\boldsymbol \vartheta}$ in the decay rate of $p_n({\boldsymbol a}).$ Mimicking the reasoning we used in the single-node case, the number of arrivals becomes Poisson with mean
\[\lambda \int_0^t\beta\left({\rm e}^{-Ru}\,{\bs\vt}\s\right)\,{\rm d}u.\]
The arrival epochs should be drawn using the density
\[f_U^{\mathbb Q}(u)=\frac{\beta\left({\rm e}^{-Ru}\,{\bs\vt}\s\right)}{\displaystyle \int_0^t\beta\left({\rm e}^{-Rv}\,{\bs\vt}\s\right){\rm d}v}.\]
Given an arrival at time $u$, $(B^{(1)},\ldots,B^{(L)})$ should be exponentially twisted by \[\big(({\rm e}^{-Ru}\,{\bs\vt}\s)_1,\ldots,
({\rm e}^{-Ru}\,{\bs\vt}\s)_L\big).\]

\subsection{Efficiency properties of importance sampling procedure}\label{S33}
We now consider the efficiency properties of the change of measure proposed in the previous subsection. 
To this end, we first argue that the vector ${\bs\vt}$ generally has some (at least one) strictly positive entries, whereas the other  entries equal 0; i.e., there are {\it no} negative entries. To this end, we first denote by ${\bs b}\s$ the `most likely point' in $A$:
\[{\bs b}\s:=\arg \inf_{{\bs b}\in A}I({\bs b}),\]
so that ${\bs \vt}\s ={\bs \vt}({\bs b}).$ It is a standard result from convex optimization that
\begin{equation}\label{DER}\frac{\partial I({\bs b})}{\partial b_i} = \vt_i({\bs b}).\end{equation}
Suppose now that $\vt_i({\bs b}\s)<0$. Increasing the $i$-th component of the ${\bs b}\s$ (while leaving all other components unchanged) would lead to a vector that is still in $A$, but that by virtue of (\ref{DER}) corresponds to a lower value of the objective function $I(\cdot)$, thus yielding that ${\bs b}$ was not the optimizer; we have thus found a contradiction. Similarly, when $\vt_i({\bs b}\s)=0$ we have that $b\s_i> a_i$ (as otherwise a reduction of the objective function value would be possible, which contradicts ${\bs b}\s$ being minimizer).

Now define $\Theta$ as the subset of $i\in\{1,\ldots,L\}$ such that $\vartheta_i>0$, and let $D\in\{1,\ldots,L\}$ the number of elements of $\Theta.$ We now argue that the number of runs needed to obtain an estimate of predefined precision scales as $n^{D/2}.$
Relying on the results from \cite{CS} (in particular their Thm.~3.4), it follows that $p_n({\bs a})$ behaves (for $n$ large) proportionally to $n^{-D/2} \exp(-nI({\bs b}))$; using the same machinery, ${\mathbb E}_{\mathbb Q}(L^2 I)$ behaves proportionally to $n^{-D/2} \exp(-2n I({\bs b}))$. Mimicking the line of reasoning of Section \ref{EFF}, we conclude that the number of runs needed  is essentially proportional to $n^{D/2}$. The formal statement is as follows. 

\begin{proposition}
As $n\to\infty$,
\begin{equation}\label{alpha}\Sigma_n\sim \alpha\,{n}^{D/2},\:\:\:\:\alpha=\frac{T^2}{\varepsilon^2} \left( \prod_{i\in D}{\vt\s_i}\right)\cdot\frac{1}{2^D}\,\left(\sqrt{2\pi}\right)^D \sqrt{\tau},\end{equation}
where $\tau$ is the determinant of the Hessian of $\log M({\bs \vt})$ in ${\bs \vt}\s.$
\end{proposition}

We further illustrate the ideas and intuition behind the qualitative result described in the above proposition by considering the case $L=2.$
It is noted that three cases may arise: (i)~$\Theta=\{1,2\}$, (ii)~$\Theta=\{1\}$, (iii)~$\Theta=\{2\}$; as case (iii) can be dealt with in the same way as case (ii), we concentrate on the cases (i) and (ii) only. In case (i), where $D=2$, the necessary condition \cite[Eqn. (3.4)]{CS}  is fulfilled as ${\bs \vt}>0$ componentwise. As in addition the conditions A--C of \cite{CS} are in place, it is concluded that \cite[Thm.\ 3.4]{CS} can be applied, leading to ${\bs b}\s = {\bs a},$ and
\[p_n({\bs a}) \sim \frac{1}{n}\frac{1}{\vt_1\s\vt_2\s\cdot 2\pi\sqrt{ \tau}} {\rm e}^{-nI({\bs a})},\]
where $\tau$ is the determinant of the Hessian of $\log M({\bs \vt})$ in ${\bs \vt}\s.$ Along the same lines, it can be shown that 
\[{\mathbb E}_{\mathbb Q}(L^2I) \sim  \frac{1}{n}\frac{1}{4\vt_1\s\vt_2\s\cdot 2\pi\sqrt{ \tau}} {\rm e}^{-2nI({\bs a})}.\]
It now follows that $\Sigma_n$ is roughly linear in $n$: with $\varepsilon$ and $T$ as introduced in Section \ref{EFF},
\begin{equation}
\label{AT1}\Sigma_n = \alpha\,n,\:\:\:\:\alpha:=\frac{T^2}{\varepsilon^2}\vt_1\s\vt_2\s\cdot\frac{\pi \,\sqrt{\tau}}{2 }.\end{equation}
In case (ii), we do not have that ${\bs \vt}>0$ componentwise, and hence \cite[Thm.\ 3.4]{CS} does not apply; in the above terminology, $D=1<2=L.$ 
Observe that in this case the exponential decay rate of the event
$\{Y^{(1)}_n(t)\geqslant na_1, Y^{(2)}_n(t)< na_2\}$
strictly majorizes that of $\{Y^{(1)}_n(t)\geqslant na_1\}$ (informally: the former event is substantially less likely than the latter).
It thus follows that $b_1\s=a_1$ and $b_2\s>a_2$, and 
\begin{eqnarray*}p_n({\bs a})&=& {\mathbb P}\left(Y^{(1)}_n(t)\geqslant na_1\right) -
 {\mathbb P}\left(Y^{(1)}_n(t)\geqslant na_1, Y^{(2)}_n(t)< na_2\right) \\
&\sim& {\mathbb P}\left(Y^{(1)}_n(t)\geqslant na_1\right) \sim \frac{1}{\sqrt{n}}\frac{1}{\vt_1\s\sqrt{2\pi\tau}}
{\rm e}^{-2nI({\bs b}\s)},
\:\:\:\tau:= \left.\frac{\rm d}{{\rm d}\vt^2} \log M(\vt,0)\right|_{\vt=\vt_1\s},\end{eqnarray*}
and in addition 
\[{\mathbb E}_{\mathbb Q}(L^2I)  \sim \frac{1}{\sqrt{n}} \,\frac{1}{2\vt_1\s\sqrt{2\pi\tau}} {\rm e}^{-2nI({\bs b}\s)}.\]
As a consequence in this regime
$\Sigma_n$ grows essentially proportional to $\sqrt{n}$ for $n$ large:
\[\Sigma_n\sim \alpha\,\sqrt{n},\:\:\:\:\alpha:=\frac{T^2}{\varepsilon^2} \, {\vt_1\s}\cdot\frac{1}{2}\,\sqrt{2\pi \tau}.\]
In case (iii) $\Sigma_n$ behaves proportionally to $\sqrt{n}$ as well. 

\subsection{Simulation experiments}\label{S34}
We conclude this section by providing a few numerical illustrations. In the first set we focus on the downstream queue only (i.e., we analyze $p_n(0,a_2)$), whereas in the second set we consider the joint exceedance probability $p_n({\bs a}).$ The precision and confidence have been chosen as in Example \ref{Ex1}. Throughout we take $t=1$, $r_1=2$, $r_2=1$, $\lambda=1$, and $\mu=1$.

\begin{figure}
\centering
%
%
%
%
%
\includegraphics[width=.9\textwidth]{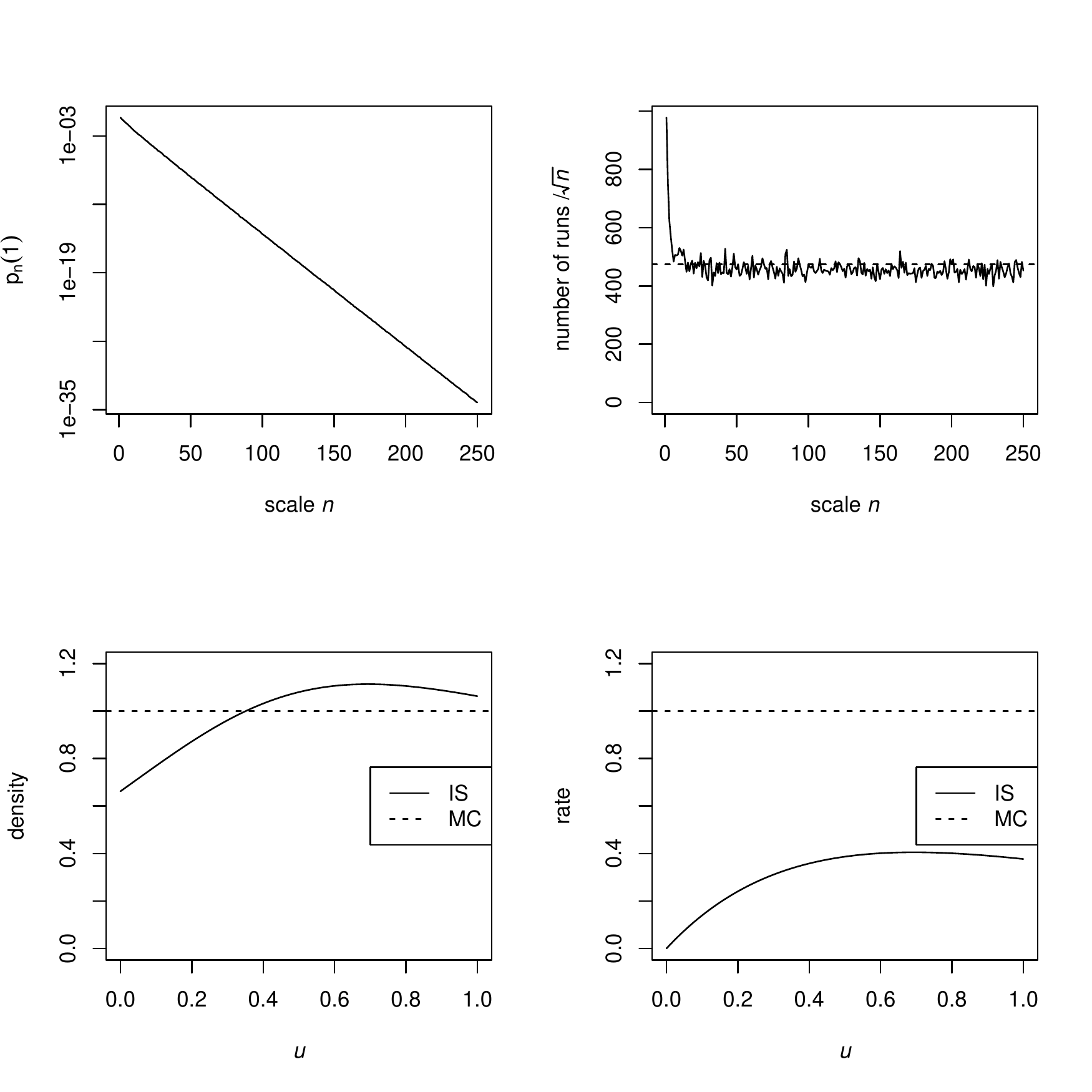}
\caption{Numerical results for Section \ref{S34}: downstream queue only.}
\label{F2}
\end{figure}

In the first set of experiments we take $a_1=0$ and $a_2=1$. Elementary numerical analysis yields that $\vt\s =0.8104$ and $\tau =1.4774$, leading to $\alpha$, as defined in (\ref{AT1}), equalling 474.3. The four panels of Fig.~\ref{F2} should be interpreted as their counterparts in Fig.\ \ref{F1}.  Interestingly, the bottom-left panel indicates that in the tandem system it does not pay off to let jobs arrive right before $t$ (as they first have to go through the first resource to end up in the second resource), as reflected by the shape of the density of the arrival epochs under ${\mathbb Q}$; to this end, recall that we reversed time in (\ref{distr}), so that a low density at $u=0$ in the graph corresponds to a high density at $u=t$ in the actual system.
The arrival rate under ${\mathbb Q}$ is $1.5103$.

\begin{figure}
\centering
%
%
%
%
%
\includegraphics[width=.9\textwidth]{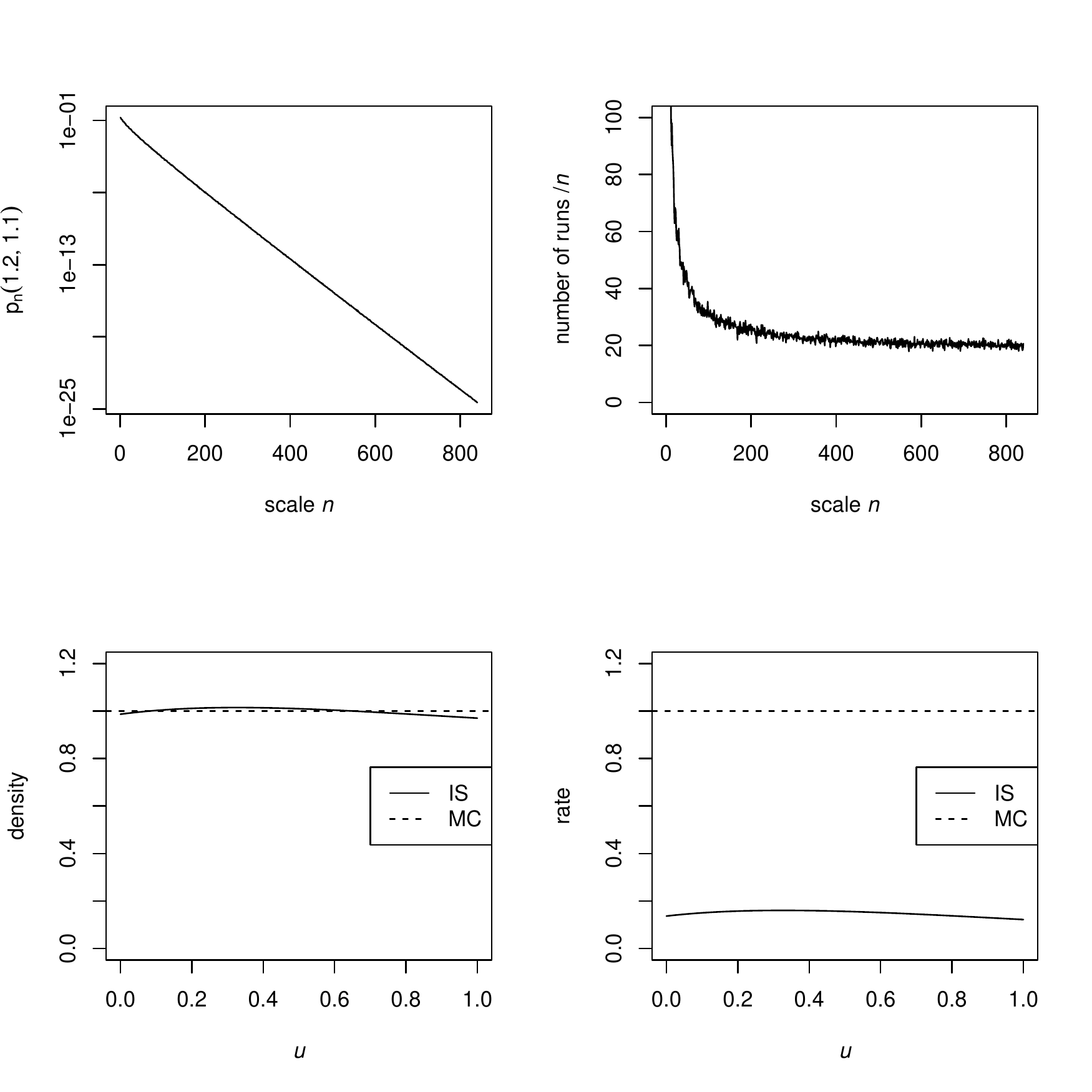}
\caption{Numerical results for Section \ref{S34}: both queues.}
\label{F3}
\end{figure}

In the second set of experiments we take $a_1=1.2$ and $a_2=1.1$; all other parameters are the same as in the first set. As mentioned above, we now consider the joint exceedance probability. As it turns out, $\vt_1\s=0.1367$ and $\vt\s_2 = 0.2225$. The top-right panel  of Fig.~\ref{F3} shows that for this specific parameter setting $\Sigma_n/n$ converges to the limiting constant rather slowly. Concerning the bottom-left panel, note that in Section \ref{S2} we saw that to make sure the first queue gets large it helps to have arrivals at the end of the interval, whereas above we observed that to make the second queue large arrivals should occur relatively early. We now focus on the event that {\it both} queues are large, and consequently the arrival distribution becomes relatively uniform again, as shown in the bottom-left panel. The arrival rate under ${\mathbb Q}$ is $2.3478$.

\section{Multi-node systems under Markov modulation}\label{S4}
In this section consider the networks analyzed in the previous section, but now in a random environment. More specifically, the type of random environment we focus on here is known as {\it Markov modulation}: the system dynamics are affected by the state of an external finite-state irreducible Markov process $J(\cdot)$ with generator matrix $Q=(q_{jj'})_{j,j'=1}^d$. When this Markov process (usually referred to as the {\it background process}) is in state $j$, arrivals occur according to a Poisson process with rate $\lambda_j$, the {\sc mgf}\:of the job size is $\beta_j({\boldsymbol \vt})$, and the routing matrix is $R_j$. Analogously to the definitions used in the case without Markov modulation, this routing matrix' $(i,i)$-th entry is \[(R_j)_{ii}:=r^{(j)}_{ii}+\sum_{i'\not=i}r^{(j)}_{ii'},\] 
which can be interpreted as the rate at which fluid leaves  server $i$ when the background process is in $j$. Likewise, for $i\not=i'$,
\[(R_j)_{ii'}:=-r^{(j)}_{ii'},\] 
which is the rate at which fluid flows from server $i$ to $i'$ when the background process is in $j$.

Below we assume that $J(0)=j_0$ for a fixed state $j_0\in\{1,\ldots,d\}$; it is seen that all results generalize to an arbitrary initial distribution in a straightforward manner.

The structure of the section is in line with that of the previous two sections: we subsequently describe general results for the model under consideration (extending earlier results from \cite{KS}), propose an importance sampling measure, establish efficiency properties of the corresponding estimator, and present a number of numerical experiments. 

\subsection{\RED{Exact expressions for moments}}\label{S41}
\RED{Before focusing on simulation-based techniques, this subsection (which can be read independently of the rest of the section) shows  that various 
moment-related quantities can be computed in closed form.}

Multi-node systems under Markov modulation have been studied in detail by Kella and Stadje \cite{KS}.
We start this \RED{subsection} by providing a compact derivation of a {\sc pde} characterizing the system's transient behavior, which was not included in that paper. To this end, we define, for $j\in\{1,\ldots,d\}$,
\[\Xi_j({\bs \vt},t) := {\mathbb E} \left(\exp\left(\sum_{\ell=1}^L \vt_\ell X^{(\ell)}(t)\right) 1_j(t)\right),\]
with $1_j(t)$ the indicator function of the event that $J(t) =j$.
Using the standard `Markov machinery', $\Xi_j({\bs \vt},t+\Delta t)$ equals (up to $o(\Delta t)$ terms) the sum of
a contribution 
\[\lambda_j \,\Delta t \,\Xi_j({\bs \vt},t) \beta_j({\bs\vt})\]
due to the scenario that an arrival occurs between $t$ and $t+\Delta t$,
a contribution 
\[\sum_{j'\not = j} q_{j'j}\, \Delta t \,\Xi_{j'}({\bs \vt},t)\]
due to the scenario  that the modulating Markov process jumps between $t$ and $t+\Delta t$,
and a contribution 
\[\left(1- \lambda_j\,\Delta t -q_{j}\,\Delta t\right)
 {\mathbb E} \left(\exp\left(\sum_{\ell=1}^L  \left(\vt_\ell
-\sum_{\ell'=1}^L\vt_{\ell'}(R_j)_{\ell\ell'} \, \Delta t 
 \right) X^{(\ell)}(t)\right)1_j(t)\right),\]
 with $q_j:=-q_{jj}$;
regarding the last term, observe that when the background process is in state $j$, and no new job arrives between $t$ and $t+\Delta t$,
\[X^{(\ell)}(t+\Delta t) = X^{(\ell)}(t)  - (R_j)_{\ell\ell} \, \Delta t \,X^{(\ell)}(t)-\sum_{\ell'\not=\ell}(R_j)_{\ell'\ell} \, \Delta t \,
X^{(\ell')}(t) .\]
We thus find that
\begin{eqnarray*}
\Xi_j({\bs\vt},t+\Delta t)&=&\lambda_j \,\Delta t \, \beta_j({\bs\vt})\,\Xi_j({\bs \vt},t)+\sum_{j'\not = j} q_{j'j}\, \Delta t \,\Xi_{j'}({\bs \vt},t)+\\
&& \left(1- \lambda_j\,\Delta t -q_{j}\,\Delta t\right)\,
\Xi_j\left({\bs \vt}- R_j{\bs\vt}\,\Delta t,t\right)+o(\Delta t).\end{eqnarray*}
This immediately leads to (by \ROOD{subsequently subtracting $\Xi_j({\bs \vt},t)$ from both sides,
dividing by $\Delta t$, and}
letting $\Delta t\downarrow 0$)
\begin{equation}\label{DV}
\frac{\partial}{\partial t} \Xi_j({\bs\vt},t) =\lambda_j \big(\beta_j({\bs\vt})-1\big)\Xi_j({\bs \vt},t)+
\sum_{j'=1}^d q_{j'j}  \,\Xi_{j'}({\bs \vt},t)-
\sum_{\ell'=1}^L \big(R_j\,{\bs \vt}\big)_{\ell'}\,\frac{\partial}{\partial \vt_{\ell'}} \Xi_j({\bs\vt},t). 
\end{equation}
Let us now compactly summarize the relation (\ref{DV}), in vector/matrix notation. This notation will prove practical when computing higher moments; in other (but related) contexts, similar  procedures have been proposed in e.g.\ \cite{GANG,RAB}.
Let ${\mathscr M}^{n_1\times n_2}$ be the set of ${\mathbb R}$-valued matrices of dimension $n_1\times n_2$ (for generic $n_1,n_2\in{\mathbb N}$). In addition, $I_n$ is the identity matrix of dimension $n\in{\mathbb N}$.
We introduce the following three matrices in ${\mathscr M}^{d\times d}$,  ${\mathscr M}^{d\times d}$, and ${\mathscr M}^{Ld\times Ld}$, respectively:
\[{\Lambda} :=\left(\begin{array}{ccc}
\lambda_1&\dots& 0\\
\vdots&\ddots& \vdots\\
0&\dots& \lambda_d\end{array}
\right),\:\:\:{B({\bs\vt})} :=\left(\begin{array}{ccc}
\beta_1({\bs\vt})&\dots& 0\\
\vdots&\ddots& \vdots\\
0&\dots& \beta_d({\bs\vt})\end{array}
\right),\:\:\:{R} :=\left(\begin{array}{ccc}
R_1&\dots& 0\\
\vdots&\ddots& \vdots\\
0&\dots& R_d\end{array}
\right).\]
We use the conventional notation $\otimes$ for the Kronecker product. Recall that the Kronecker product is bilinear, associative and distributive with respect to addition; these properties we will use in the sequel without mentioning. It also satisfies the mixed product property $(A\otimes B)(C\otimes D)=(AC)\otimes(BD)$. Furthermore, note that $I_{n_1}\otimes I_{n_2} = I_{n_1n_2}$. 

We now consider some differentiation rules for matrix-valued functions which will allow us to iteratively evaluate moments. In the first place we
define the operator  $\nabla_{\bs \vt}$ for ${\bs\vt}\in{\mathbb R}^L$; to keep notation compact, we often suppress the subscript ${\bs\vt}$, and write just $\nabla$. Let $f\equiv f({\bs \vt})$ be  a mapping of ${\mathbb R}^L$ to 
$\mathscr{M}^{n_1\times n_2}$. Then $\nabla f\equiv \nabla f({\bs\vt})\in{\mathscr M}^{n_1L\times n_2}$ is defined by 
\[
\nabla f=\begin{pmatrix}
\nabla f_{11} & \nabla f_{12} &\cdots &\nabla f_{1n_2}\\
\nabla f_{21} & \nabla f_{22} &\cdots &\nabla f_{2n_2}\\
\vdots             & \vdots & \ddots & \vdots\\
\nabla f_{n_11} & \nabla f_{n_12} & \cdots & \nabla f_{n_1n_2}
\end{pmatrix},\:\:\:
\mbox{where}\:\:\:
\nabla f_{ij} := 
\begin{pmatrix}
\partial_1 f_{ij}\\
\partial_2 f_{ij}\\
\vdots\\
\partial_L f_{ij}
\end{pmatrix}.
\]
In the above definition  $\nabla f_{ij}\equiv \nabla f_{ij}({\bs\vt})$ is to be understood as the usual gradient; the symbol $\partial_i$ is used to denote the partial derivative with respect to the $i$-th variable, in the sense
of
\[\partial_i f_{ij}:= \frac{\partial}{\partial \vt_i}f_{ij}(\bs\vt).\]
Furthermore, we define inductively $\nabla^k f\equiv \nabla^k f({\bs\vt}) := \nabla(\nabla^{k-1} f)$, $k\in{\mathbb N}$, with $\nabla^0 f := f$. It is checked that $\nabla^k f({\bs\vt})$ is a mapping of ${\mathbb R}^L$ to 
$\mathscr{M}^{L^kn_1\times n_2}$. 

In the sequel we use a couple of differentiation rules, that we have listed below. Let $A(\cdot)$ be a matrix-valued function from ${\mathbb R}^L$ to ${\mathscr M}^{n_1\times n_2}$, and $B(\cdot)$ a matrix-valued function from ${\mathbb R}^L$ to ${\mathscr M}^{n_2\times n_3}$, and let $I_q$ be a $q\times q$ identity matrix  (for some $q\in\mathbb{N}$). Then,
\begin{itemize}
\item[--] {\it Product rule}:
\[ \nabla_{\bs\vt} \big(A({\bs\vt}) \,B({\bs \vt})\big) = (\nabla_{\bs\vt} A({\bs\vt}))\, B({\bs\vt}) + (A(\bs \vt)\otimes I_L) \, \nabla_{\bs\vt} B({\bs\vt});\]
being an element of $ {\mathscr M}^{L n_1\times n_3}$.
\item[--] {\it  Differentiation of Kronecker product (1)}:
\[\nabla_{\bs\vt}(I_q \otimes A(\bs\vt)) =I_q \otimes (\nabla_{\bs\vt} A(\bs\vt)).\]
\item[--] {\it Differentiation of Kronecker product (2)}:
\begin{eqnarray*}\nabla_{\bs\vt}(A(\bs\vt)\otimes I_q) &=& (K_{n_1,q}\otimes I_L)(I_q\otimes(\nabla_{\bs\vt} A(\bs\vt)))K_{q,n_2} \\&=& (K_{n_1,q}\otimes I_L)K_{q,n_2}(\nabla_{\bs\vt}A(\bs\vt)\otimes I_q),\end{eqnarray*} where $K_{m,n}$ is the commutation matrix defined by 
\[
K_{m,n} :=\sum_{i=1}^m\sum_{j=1}^n (H_{ij} \otimes H^{\rm T}_{ij}),
\]
and $H_{ij}\in\mathscr{M}^{m\times n}$ denotes a matrix with a $1$ at its $(i,j)$-th position and zeros elsewhere,  cf.\ \cite{MN1979}. 
\end{itemize}
The first rule can be checked componentwise and the second rule is trivial. The third rule follows from the first and second rule in combination with the fact that the Kronecker product commutes after a correction with the commutation matrices. Moreover, we use the property $K_{m,n}^{-1}=K_{n,m}$. An overview of the properties of commutation matrices can be found in \cite{MN1979}. 

In the introduced terminology, it follows that (\ref{DV}) can be written as 
\begin{equation}\label{partdifeq}\frac{\partial}{\partial t} {\bs \Xi}({\bs\vt},t) = \Lambda \big(B({\bs\vt}) -I_d\big)\,{\bs \Xi}({\bs\vt},t) +
Q^{\rm T}\, {\bs \Xi}({\bs\vt},t) -
\big(I_d\otimes {\bs \vt}^{\rm T}\big) R^{\rm T}\, \nabla_{\bs\vt}{\bs \Xi}({\bs\vt},t) .\end{equation}

We now point out how (transient and stationary) moments can be evaluated; note that \cite{KS} focuses on the first two stationary moments at epochs that the background process jumps. 
We throughout use the notation ${\bs z}_i(t)$ for the $i$-th derivative of ${\bs \Xi}({\bs\vt},t)$ in $({\bs 0},t)$, for $t\geqslant 0$:
\[{\bs z}_i(t) := \left.\nabla^i_{\bs\vt}{\bs\Xi}({\bs\vt},t)\right|_{{\bs\vt}=\bs0}\in{\mathscr M}^{L^id\times d},\] 
for $i\in{\mathbb N}$. Note that, with $\pi_j(t)= (\exp(Qt))_{j_0,j}$,
\[\bs\Xi({\bs\vt},0)={\bs e}_{j_0},\:\:\:\bs\Xi({\bs0},t) = {\bs\pi}(t)^{\rm T}\equiv (\pi_1(t),\ldots,\pi_d(t)).\]

\newcommand{\f}[2]{\frac{#1}{#2}}

\vb

$\circ$\:\:
We start by characterizing the first moments. 
Applying the operator $\nabla \equiv \nabla_{\bs\vt}$  to the differential equation (\ref{partdifeq}) yields
\begin{eqnarray}\nonumber
\nabla_{\bs\vt}\left(\f{\p}{\p t} {\bs\Xi}({\bs\vt},t)\right)&=&(\Lambda\otimes I_L)(\nabla_{\bs\vt} B(\bs\vt))\,{\bs\Xi}(\bs\vt,t)\,+
\\
\nonumber &&\left(Q^{\rm T}\otimes I_L +\Lambda(B(\bs\vt)-I_d)\otimes I_L - R^{\rm T}\right)\nabla_{\bs\vt}{\bs\Xi}(\bs\vt,t)\,-\\
\label{z1}&&\left(((I_d\otimes{\bs\vt}^{\rm T})R^{\rm T}) \otimes I_L\right)\, \nabla_{\bs\vt}^2\,{\bs\Xi}(\bs\vt,t),
\end{eqnarray}
using standard properties of the Kronecker product in combination with
\[
\nabla_{\bs\vt}(I_d\otimes {\bs\vt}^{\rm T})= I_d\otimes(\nabla_{\bs\vt}{\bs\vt}^{\rm T}) = I_d\otimes(\bs  e_1,\ldots,\bs  e_L) = I_d\otimes I_L =  I_{dL},
\]
where $\bs e_i$ denotes the $L$-dimensional column vector in which component $i$ equals $1$ and all other components are $0$. 
Then, inserting ${\bs\vt}=\bs0$ into (\ref{z1}) yields the system of (non-homogeneous) linear differential equations
\begin{equation}
\label{eq:z1}
\bs z_1'(t) = (\Lambda\otimes I_L)\,\nabla B(\bs0)\,\bs \pi(t) + \big((Q^{\rm T} \otimes I_L)-R^{\rm T}\big)\,\bs z_1(t).
\end{equation}
In the stationary case, we obtain
\begin{equation}
\label{eq:z1stat}
\bs z_1(\infty) = \big(R^{\rm T}-(Q^{\rm T} \otimes I_L)\big)^{-1}(\Lambda\otimes I_L)\,\nabla B(\bs0)\,\bs \pi,
\end{equation}
with $\bs\pi := \lim_{t\to\infty}\bs \pi(t)$ (being the unique non-negative solution of ${\bs\pi}^{\rm T}Q={\bs 0}^{\rm T}$ such that the entries of ${\bs \pi}$ sum to $1$).

\vb

$\circ$\:\:
We now move to second moments. Applying the $\nabla_{\bs\vt}$-operator to (\ref{z1}), 
\begin{eqnarray*}
\nabla^2_{\bs\vt}\left(\f{\p}{\p t}\bs\Xi(\bs\vt,t)\right)&=&(\Lambda\otimes I_{L^2})(\nabla_{\bs\vt}^2 B(\bs\vt))\bs\Xi(\bs\vt,t)+\\&& \big(((\Lambda\otimes I_L)\nabla_{\bs\vt} B(\bs\vt))\otimes I_L\big)\nabla_{\bs\vt}\bs\Xi(\bs\vt,t)+\\&&
\nabla_{\bs\vt}(\Lambda B(\bs\vt)\otimes I_L)\nabla_{\bs\vt}\bs\Xi(\bs\vt,t)+\\
&&\big(Q^{\rm T}\otimes I_{L^2}+\Lambda(B(\bs\vt)-I_d)\otimes I_{L^2}- R^{\rm T}\otimes I_L\big)\,\nabla^2_{\bs\vt}\bs\Xi(\bs\vt,t)-\\
&&(((I_d\otimes\bs\vt^{\rm T})R^{\rm T})\otimes I_{L^2})\,\nabla^3_{\bs\vt}\bs\Xi(\bs\vt,t)-\\&&\nabla_{\bs\vt}(((I_d\otimes\bs\vt^{\rm T})R^{\rm T})\otimes I_{L})\,\nabla^2_{\bs\vt}\bs\Xi(\bs\vt,t),
\end{eqnarray*}
in which the factor $\nabla_{\bs\vt}(\Lambda B(\bs\vt)\otimes I_L)$ can be expressed more explicitly as
\[
(K_{d,L}\otimes I_L)K_{L,dL}(((\Lambda\otimes I_L)\nabla_{\bs\vt}B(\bs\vt))\otimes I_L),
\]
and the factor $\nabla_{\bs\vt}(((I_d\otimes\bs\vt^{\rm T})R^{\rm T})\otimes I_{L})$ simplifies to
$
(K_{d,L}\otimes I_L)K_{L,dL}(R^{\rm T}\otimes I_L)
$.
Inserting $\bs\vt=0$ yields the system of linear differential equations
\begin{eqnarray*}
\bs z_2'(t) &=& (\Lambda\otimes I_{L^2})\,(\nabla^2 B(\bs 0))\,\bs \pi(t)+\\&&
(Q^{\rm T}\otimes I_{L^2}-((K_{d,L}\otimes I_L)K_{L,dL}+I_{dL^2})(R^{\rm T}\otimes I_L)) \,\bs z_2(t)+\\
&&\big(((\Lambda\otimes I_L)(\nabla B(\bs 0)))\otimes I_L\big)\bs z_1(t) +\\&& (K_{d,L}\otimes I_L)K_{L,dL}(((\Lambda\otimes I_L)\nabla B(\bs0))\otimes I_L) \bs z_1(t)
\end{eqnarray*}
where $\bs z_1(t)$ solves Eqn.\ \eqref{eq:z1}. As before, the stationary quantities can be easily derived (by equating $\bs z_2'(t)$ to $0$). One has to keep in mind, however, that some of the mixed partial derivatives occur multiple times in $\bs z_k$, for $k\in\{2,3,\ldots\}$, and therefore the solution will only be unique after removing the corresponding redundant rows. Alternatively, the system can be completed by including equations which state that these mixed partial derivatives are equal. 

\vb

$\circ$\:\:
It now follows that higher moments can be found recursively,  using the three differentiation rules that we stated above.

\begin{remark}{\em Various variants of our model can be dealt with similarly. 
In this remark we consider the slightly adapted model in which shots only occur simultaneously with a jump in the modulating Markov chain. Then (up to $o(\D t)$ terms) $\Xi_j(\bs\vartheta,t+\D t)$ is the sum of a contribution
\[
\sum_{j'\neq j} q_{j'j} \D t \,\Xi_{j'}(\bs\vartheta,t) \beta_j(\bs\vartheta)
\]
due to the scenario that there is a jump in the modulating chain in the interval $[t,t+\Dt]$ (which also induces a shot), and a contribution of
\[
(1-q_j\D t)\,\E\bigg(\exp\Big(\sum_{\ell=1}^L\Big(\vartheta_\ell-\sum_{\ell'=1}^d \vartheta_{\ell'}(R_j)_{\ell\ell'}\D t\Big)X^{(\ell)}(t)\Big)1_{j}(t)\bigg),
\]
with $q_j:=-q_{jj}$, in the scenario that there is no jump. 
Performing the same steps as above, we obtain
\[
\f{\p}{\p t} \Xi_j({\bs \vt},t) =q_{j}(\beta_j({\bs \vt})-1)\Xi_j({\bs \vt},t)+ \sum_{j'=1}^d q_{j'j} \Xi_{j'}({\bs \vt},t) \beta_j({\bs \vt}) - \sum_{j'=1}^L (R_j{\bs \vt})_{j'} \f{\p}{\p\vartheta_{j'}}\Xi_j({\bs \vt},t),\]
which has a similar structure as (\ref{DV}).
It follows that the moments can be found as before.
With $\widetilde{Q}:={\rm diag}\{q_1,\ldots,q_d\}$, it turns out that the transient means are given by
\[
\bs z_1'(t) =  \nabla B(\bs 0)( Q^{\rm T} +\widetilde{Q})\bs\pi(t) + \big((Q^{\rm T}\otimes I_L) - R^{\rm T}\big)\bs z_1(t).
\]
In particular, the stationary first moment equals
\[
\bs z_1(\infty) = \big(R^{\rm T}-(Q^{\rm T}\otimes I_L)\big)^{-1} \nabla B(\bs 0) ( Q^{\rm T} +\widetilde{Q})\bs\pi.
\]
}
\end{remark}

Consider the following numerical example for the computation of the expected values and variances, in which the technique described above is illustrated.

\begin{figure}
\centering
\includegraphics[width=.9\textwidth]{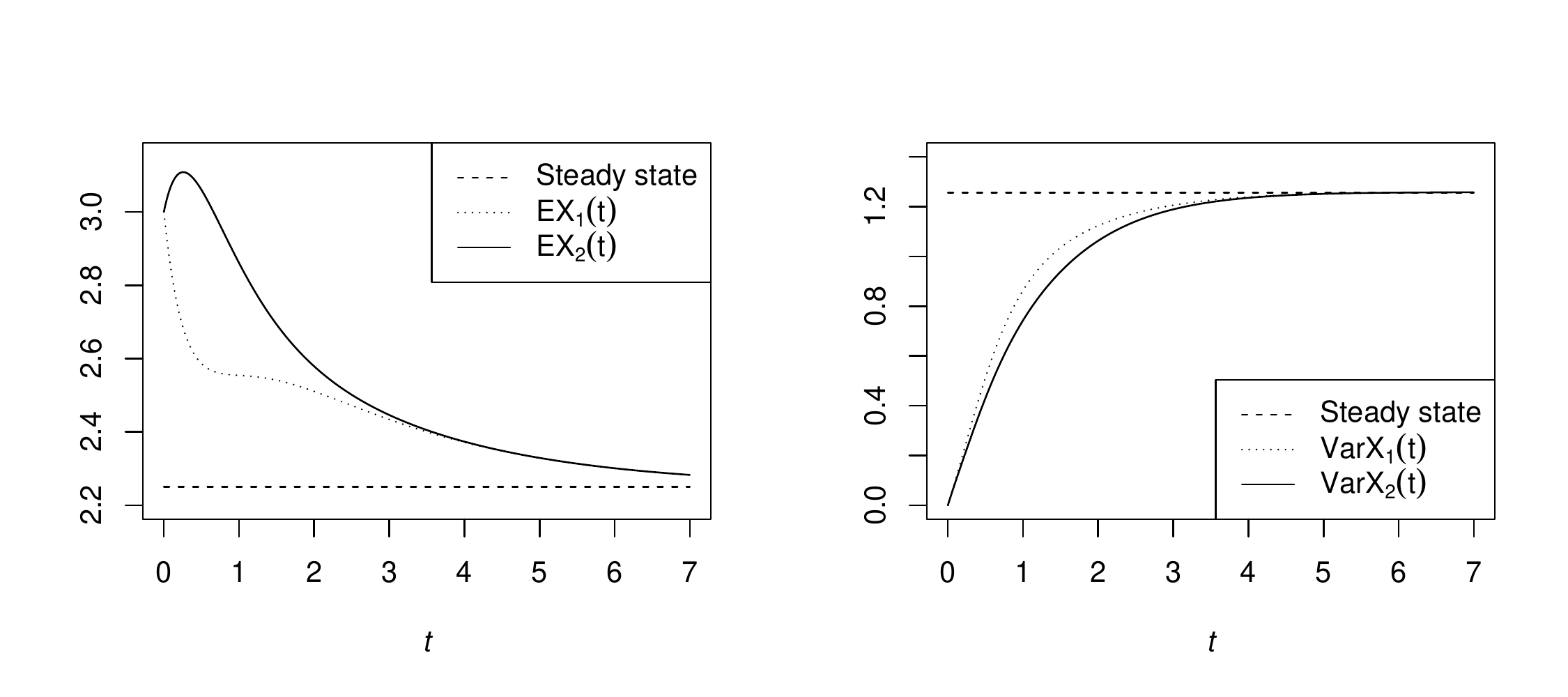}
\caption{\textit{Transient expected values and variances of Example \ref{EX2}.}}
\label{fig:EV}
\end{figure}

\begin{figure}
\centering
\includegraphics[width=.7\textwidth]{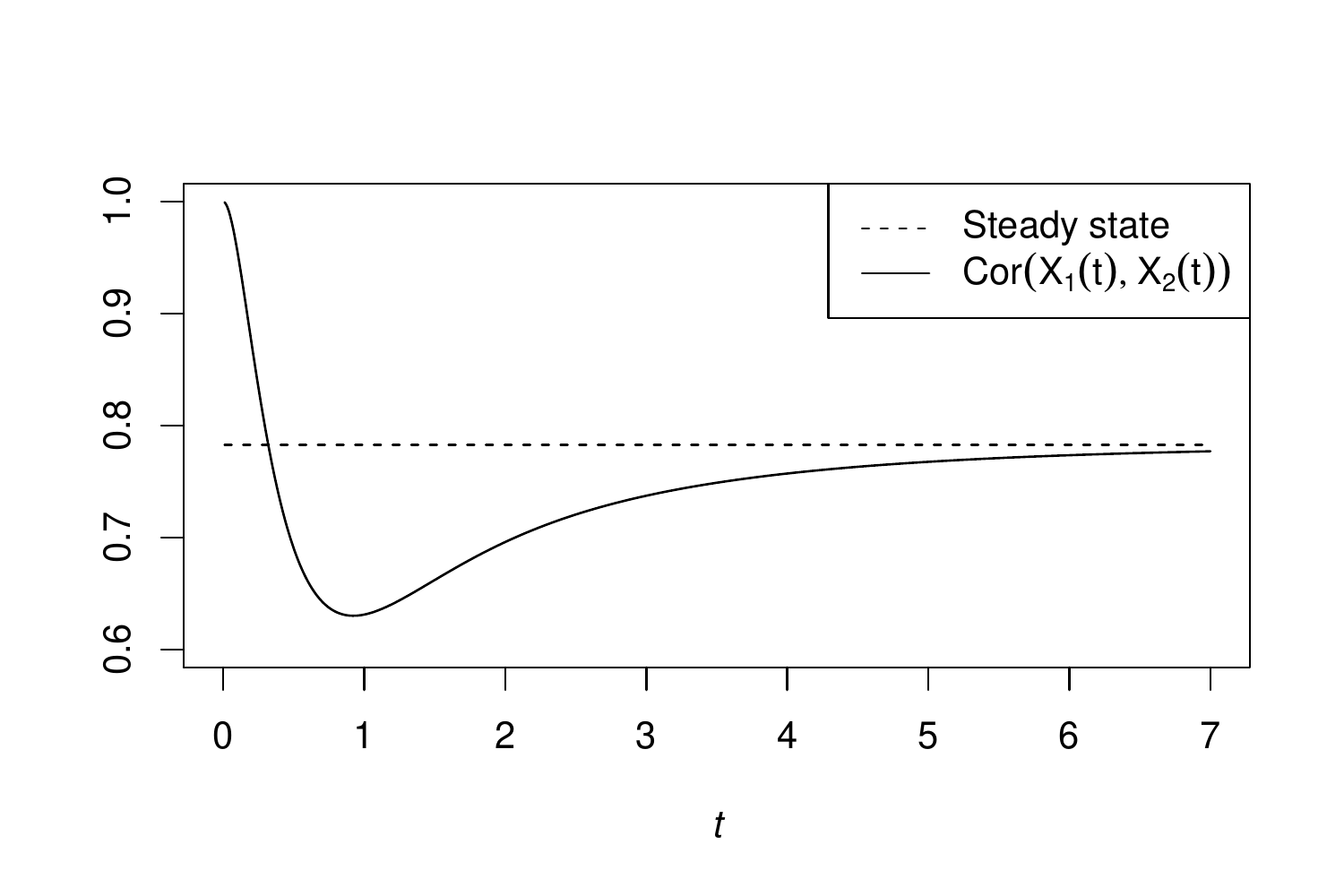}
\vspace*{-3mm}
\caption{\textit{Transient correlation between $X_1(t), X_2(t)$ of Example $\ref{EX2}$.}}
\label{fig:Corr}
\end{figure}

\begin{example}\label{EX2}
{\em 
In this example, we choose the parameters in such a way that we see non-monotonic behavior. Our example corresponds to a case in which the system does not start empty, which is dealt with by imposing suitable starting conditions. We consider a two-dimensional ($L=2$) queueing system, with a two-dimensional state space of the Markov modulating process ($d=2$). We pick $q_{12}=q_{21}=1$, $\lambda_1=\lambda_2=1$, ${\mathbb E}\, B_1={\mathbb E}\, B_2={\mathbb E}\, B_1^2={\mathbb E}\, B_2^2=1$, $J(0)=1$, $(X_1(0),X_2(0))=(3,3)$, and the rate matrices 
\[
R_1 = 
\begin{pmatrix}
2&-1\\
-1&1
\end{pmatrix}, \qquad\qquad
R_2 =
\begin{pmatrix}
1&-1\\
-1&2
\end{pmatrix}.
\]
Due to  the symmetry in the choice of the parameters, one can expect that for both states of the background process expected value tends (as $t$ grows large) to the same steady-state value; the same is anticipated for the stationary variance. This is confirmed by 
Fig.~\ref{fig:EV}. For $t$ small, the two queues behave differently due to $J(0)=1$, which implies that queue~$1$ drains faster. Note that ${\mathbb E}\, X_2(t)$ even increases for $t$ small, due to the fact that the flow from node~1 to 2 equals the flow from 2 to 1, constituting a net flow of zero, so that the additional contribution of external output to node 2 leads to a net increase of ${\mathbb E}\, X_2(t)$.  The transient correlation is plotted in Fig.\ \ref{fig:Corr}. At time $t=0$ the queues are perfectly correlated, since the starting state is deterministic. Then the correlation decreases due to the asymmetric flow rates until around $t=1$, which is when the Markov chain $J$ is expected to switch, after which the correlation monotonously tends to the steady state.}
\end{example}

\subsection{Tail probabilities, change of measure}\label{S42}
We now characterize the decay rate of the rare-event probability under study, and we propose a change of measure to efficiently estimate it. In the notation we have been using so far, we again focus on
\[p_n({\boldsymbol a}) :=
{\mathbb P}\left(Y^{(1)}_n(t)\geqslant na_1, \ldots,Y^{(L)}_n(t)\geqslant na_L\right)=
{\mathbb P}\left({\bs Y}\hspace{-0.5mm}_n(t)\in A\right),
\]
where ${\bs Y}\hspace{-0.5mm}_n(t)=(Y_n^{(1)}(t),\ldots,Y_n^{(L)}(t)).$ \RED{It is stressed that, following \cite{BM}, we consider the regime in which the background process is `slow'. In concrete terms, this means that we linearly speed up the driving Poisson process (i.e., we replace the arrival rates $\lambda_j$ by $n\lambda_j$), but leave the rates of the Markovian background process unaltered.}

\newcommand{\hf}{^{(f)}}
\newcommand{\hff}{{(f)}}

First we find an alternative characterization of the state of the system at time $t$. 
Let ${\mathscr F}_t$ denote the set of all functions from $[0,t]$ onto the states $\{1,\ldots, d\}$. Consider a path $f\in{\mathscr F}_t$. Let $f$ have $K\hff$ jumps between $0$ and $t$, whose epochs we denote by $t_1\hff$ up to $t_{K\hff}\hff$ (and in addition $t_0\hff:=0$ and $t_{K\hff+1}\hff:=t$). Let \[j_i(f):=\lim_{t\downarrow t_i(f)}f(t)\] (i.e., the state of $f$ immediately after the $i$-th jump). 
We also introduce
\[D_i(u,f):=\exp\left(-(t_{i+1}(f)-u)\,R_{j_i(f)}\right),\:\:\:
D_i(f):=\exp\left(-({t_{i+1}(f) }-{t_i(f)})\,R_{j_i(f)}\right)
.\]
Suppose now that the Markov process $J(\cdot)$ follows the path $f\in{\mathscr F}_t$. Then the contribution to the {\sc mgf}\:of ${\bs  X}(t)$ due to shots that arrived between $t_i(f)$ and $t_{i+1}(f)$ is,  mimicking the arguments that we used in Section \ref{S32} for non-modulated networks,
\[\psi_i(f,{\bs \vt}):= \exp\left(\lambda_{j_i(f)} \int_{t_i(f)}^{t_{i+1}(f)}\big( \beta_{j_i(f)}\big(D_i(u,f) D_{i+1}(f)\cdots D_{K(f)}(f)\,{\bs \vt}\big)-1\big){\rm d}u\right).\] As a consequence, the {\sc mgf}\:of $X(t)$ given the path $f$ is 
\[M_f({\bs \vt}):= \prod_{i=0}^{K(f)} \psi_i(f,{\bs \vt}).\]
First conditioning on the path of $J(\cdot)\in{\mathscr F}_t$ between $0$ and $t$ and then unconditioning, it then immediately follows that the {\sc mgf}\:of ${\bs X}(t)$ is given by
\[M({\bs \vt}) = {\mathbb E}\,M_J({\bs\vt}).\]

Then, precisely as is shown in \cite{BM} for a related stochastic system, the decay rate can be characterized as follows:
\begin{equation}
\label{drMM}
\lim_{n\to\infty} \frac{1}{n}\log p_n({\bs a}) = -\inf_{f\in {\mathscr F}_t}{\mathbb I}_f({\bs a}),\:\:\:
{\mathbb I}_f({\bs a}):= \inf_{{\bs b}\in A} \sup_{\bs \vt}\left(\langle {\bs \vt},{\bs b}\rangle - \log
M_f({\bs \vt}) \right).\end{equation}
\ROOD{The argumentation to show this is analogous to the one in \cite[Thm.\ 1]{BM}, and can be summarized as follows. In the first place, let $f\s$ be the optimizing path in (\ref{drMM}). Then, as $J(\cdot)$ does not depend on $n$, we can choose a `ball' ${\mathscr B}_t({f\s})$ around $f\s$ such that the decay rate of the probability of $J(\cdot)$ being in that ball is $0$. The lower bound follows by only taking into account the contribution due to paths in ${\mathscr B}_t({f\s})$. The upper bound follows by showing that the contribution of all $f\in {\mathscr F}_t\setminus
{\mathscr B}_t({f\s})$ is negligible.} 

Informally, the path $f\s$ has the interpretation of the most likely path of $J(\cdot)$ given that the rare event under consideration happens.
To make sure that the event under consideration is rare, we assume that for all $f\in{\mathscr F}_t$
\[\left(\left.\frac{\partial}{\partial \vt_1}M_f({\bs \vt})\right|_{{\bs\vt}={\bs 0}},\ldots,
\left.\frac{\partial}{\partial \vt_L}M_f({\bs \vt})\right|_{{\bs\vt}={\bs 0}}\right)\not \in A.\]

\vb

The change of measure we propose is the following. \ROOD{In every run we first sample} the path $J(s)$ for $s\in[0,t]$ {\it under the original measure} ${\mathbb P}$ (i.e., with $J(0)=j_0$, and then using the generator matrix $Q$). We call the resulting path $f\in{\mathscr F}_t$. 
For this path, define ${\bs\vt}\s_f\geqslant {\bf 0}$ as the optimizing $\bs\vt$ in the definition of ${\mathbb I}(f)$ in (\ref{drMM}); ${\bs b}\s_f\in A$ is the optimizing ${\bs b}.$

Conditional on the path $f$ of the background process, under the new measure ${\mathbb Q}$ the number of external arrivals  between $t_{i}(f)$ and $t_{i+1}(f)$ is Poisson with parameter
\[\int_{t_i(f)}^{t_{i+1}(f)}\lambda_{j_i(f)} \beta_{j_i(f)}  \left(P_i(u,f)\,{\bs\vt}\s_f\right){\rm d}u,\]
where $P_i(u,f):=D_i(u,f) D_{i+1}(f)\cdots D_{K(f)}(f).$
The arrival epochs between $t_{i}(f)$ and $t_{i+1}(f)$ should be drawn using the density
\[f_U^{\mathbb Q}(u)=\frac{\beta_{j_i(f)} \left(P_i(u,f)\,{\bs\vt}\s_f\right)}{\displaystyle \int_{t_i(f)}^{t_{i+1}(f)}\beta_{j_i(f)} \left(P_i(v,f)\,{\bs\vt}\s_f\right){\rm d}v}.\]
Given an arrival at time $u$ between $t_{i}(f)$ and $t_{i+1}(f)$, the job sizes $(B^{(1)},\ldots,B^{(L)})$
should be sampled from a distribution with {\sc mgf}\:$\beta_{j_i(f)}({\bs\vt})$, but then exponentially twisted by 
\[\left(\left(P_i(u,f)\,{\bs\vt}\s_f\right)_1,\ldots,
\left(P_i(u,f)\,{\bs\vt}\s_f\right)_L\,\right).\]

\ROOD{\begin{remark} {\em
As mentioned above, the background process is sampled under the original measure, whereas an alternative measure is used for the number of arrivals, the arrival epochs, and the job sizes. The intuition behind this, is that the rare event under consideration is caused by two effects:
\begin{itemize}
\item[$\circ$] In the first place, samples of the background process $J$ should be close to $f\s$. Under ${\mathbb P}$ a reasonable fraction ends up close to $f\s$ --- more precisely, the event of $J$ being close to $f\s$ does not become increasingly rare when $n$ grows. As a consequence, no change of measure is needed here. 
\item[$\circ$] In the second place, given the path of the background process, the $Y^{(\ell)}_n(t)$  should exceed the values $na_\ell$, for $\ell=1,\ldots,L$. This event {\it does} become exponentially rare as $n$ grows, so importance sampling is to be applied here. 
\end{itemize}}
\end{remark}}

\subsection{Efficiency properties of importance sampling procedure}\label{S43}

We now analyze the speed up realized by the change of measure introduced in the previous subsection. Unlike our results for the non-modulated systems, now we cannot find the precise rate of growth of $\Sigma_n.$ What {\it is} possible though, is proving {\it asymptotic efficiency} (also sometimes referred to as {\it logarithmic efficiency}), in the sense that we can show that 
\[\lim_{n\to\infty} \frac{1}{n}\log {\mathbb E}_{\mathbb Q}(L^2I) 
=\lim_{n\to\infty} \frac{2}{n}\log p_n({\bs a}) = - 2\inf_{f\in {\mathscr F}_t} \inf_{{\bs b}\in A} \sup_{\bs \vt}\left(\langle {\bs \vt},{\bs b}\rangle - \log
M_f({\bs \vt}) \right)\]
(where the second equality is a consequence of (\ref{drMM})). This equality is proven as follows. As by Jensen's inequality  ${\mathbb E}_{\mathbb Q}(L^2I) \geqslant
({\mathbb E}_{\mathbb Q}(LI))^2 =  (p_n({\bs a}))^2,$ we are left to prove the upper bound:
\[\lim_{n\to\infty} \frac{1}{n}\log {\mathbb E}_{\mathbb Q}(L^2I) 
\leqslant \lim_{n\to\infty} \frac{2}{n}\log p_n({\bs a}).\]
If the path of $J(\cdot)$ equals $f\in{\mathscr F}_t$,
it follows by an elementary computation that we have constructed the measure ${\mathbb Q}$ such that
\[L
= \frac{{\rm d}{\mathbb P}}{{\rm d}{\mathbb Q}} =\prod_{\ell =1}^L  \exp\left(-\langle{\bs \vt}\s_f, {\bs Y}\hspace{-0.5mm}_n(t)\rangle+n \,\log M_f({\bs\vt}\s_f)\right).\]
The fact that ${\bs\vt}\s_f$ is componentwise non-negative, in combination with the fact  that
${\bs Y}\hspace{-0.5mm}_n(t)\geqslant {\bs a}$ when $I=1$, entails that
\[
LI \leqslant \exp\left(-n\,\langle{\bs\vt}\s_f,{\bs a}\rangle+n \,\log M_f({\bs\vt}\s_f) \right)=
\exp\left(-n\,\langle{\bs\vt}\s_f,{\bs b}\s_f\rangle + n \,\log M_f({\bs\vt}\s_f)\right)={\rm e}^{-n\,{\mathbb I}_f({\bs a})},\]
noting that ${\bs a}$ and ${\bs b}\s_f$ may only differ if the corresponding entry of ${\bs\vt}\s_f$ equals $0$ (that is,  
$\langle {\bs a}-{\bs b}_f\s,{\bs\vt}_f\s\rangle =0$).
The upper bound thus follows: with $f\s$  the minimizing path in (\ref{drMM}), recalling that $J(\cdot)$ is sampled under ${\mathbb P}$,
\[
{\mathbb E}_{\mathbb Q}(L^2I)\leqslant {\mathbb E}\,
{\rm e}^{-2n\,{\mathbb I}_J({\bs a})}\leqslant 
{\rm e}^{-2n\,{\mathbb I}_{f\s}({\bs a})}.\]
We have established the following result.
\begin{proposition}
As $n\to\infty$, the proposed importance sampling procedure is asymptotically efficient. This means that the number of runs needed grows subexponentially:
\[\lim_{n\to\infty} \frac{1}{n}\log \Sigma_n = 0.\]
\end{proposition}

\begin{figure}
\centering
\includegraphics[width=.9\textwidth]{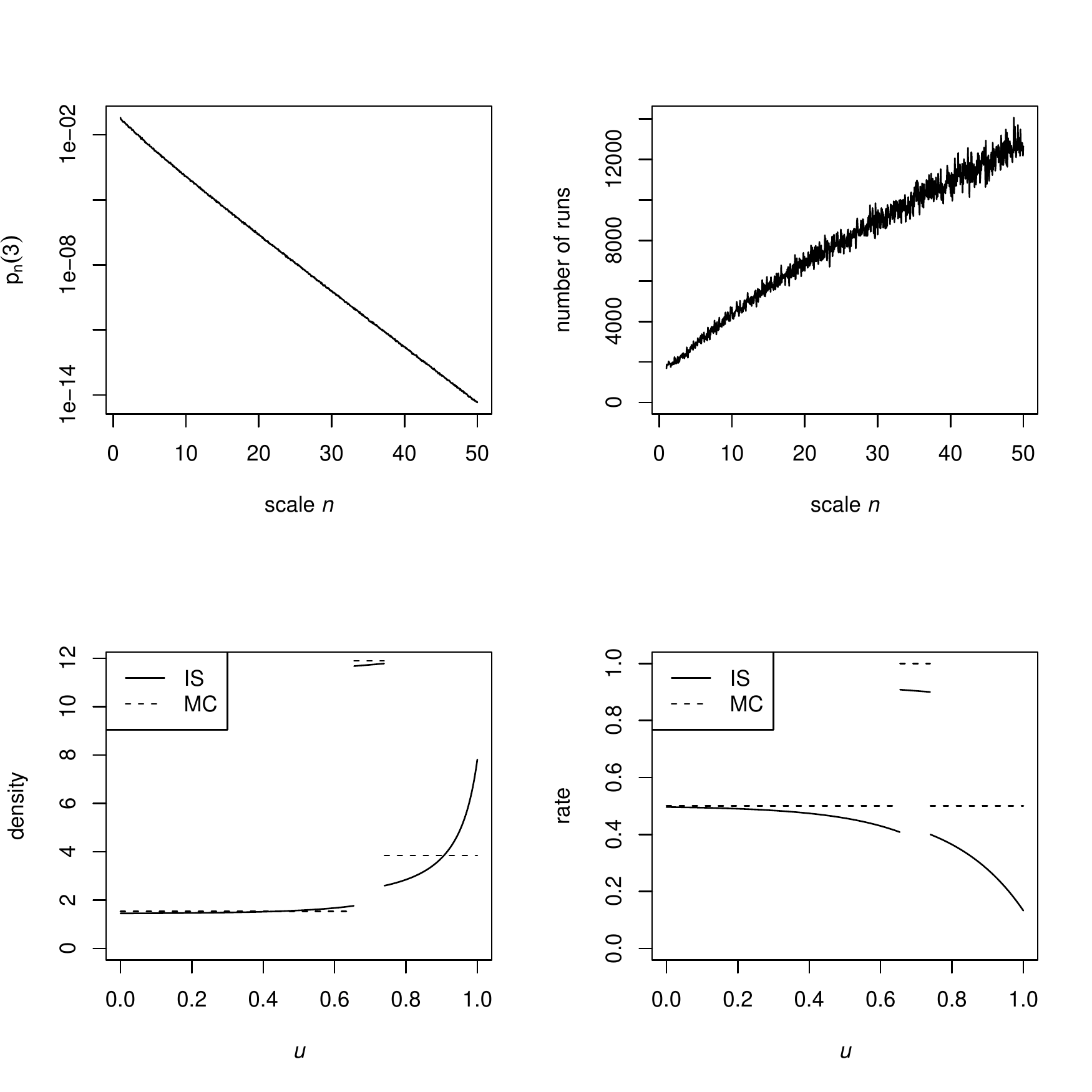}
\caption{Numerical results for Section \ref{S44}: first example.}
\label{F4}
\end{figure}

\begin{figure}
\centering
\includegraphics[width=.9\textwidth]{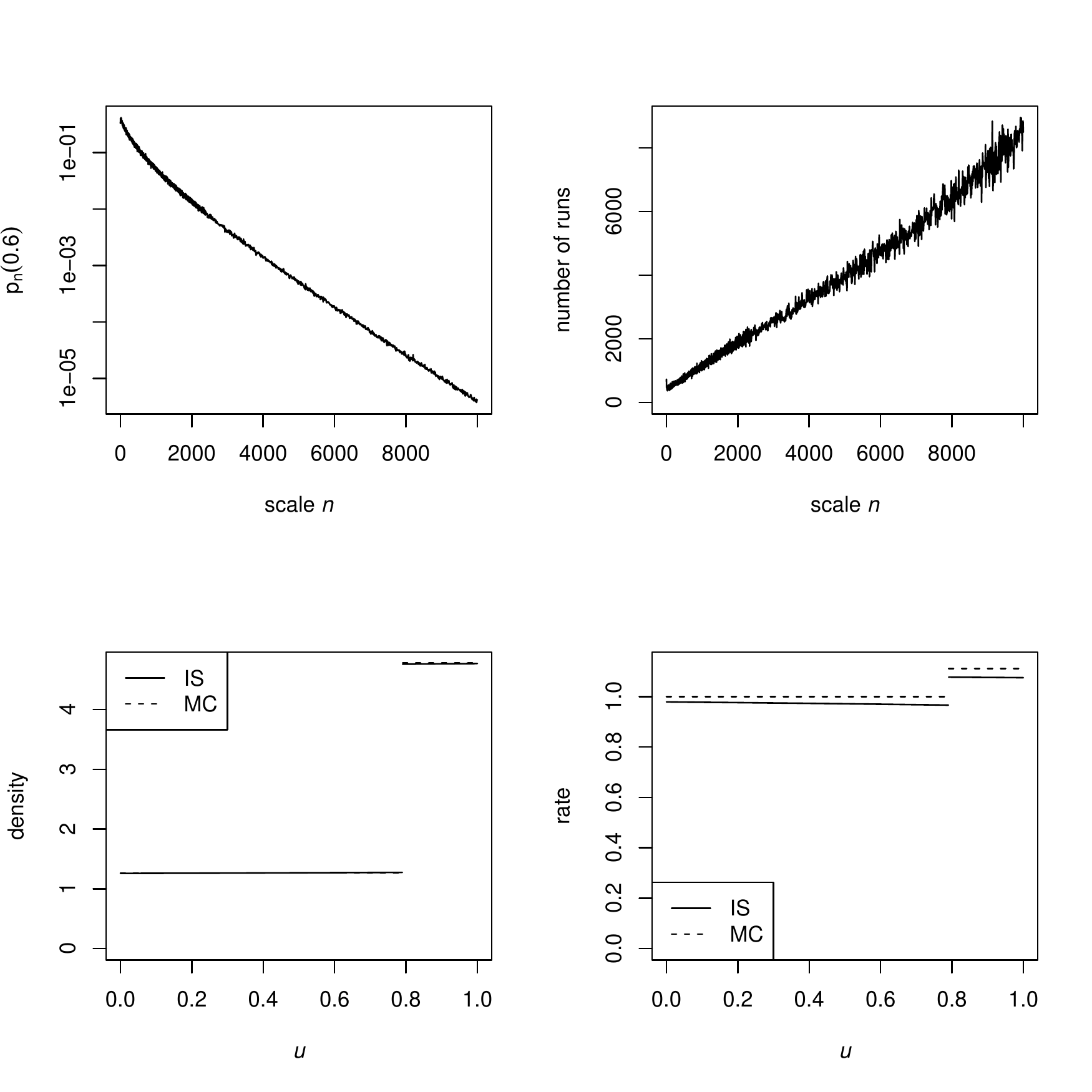}
\caption{Numerical results for Section \ref{S44}: second example.}
\label{F5}
\end{figure}

\begin{remark} \RED{\em In the scaling considered, for both the logarithmic asymptotics of $p_n({\bs a})$ and  our importance sampling algorithm,
the precise transition rates $q_{ij}$ do not matter; the only crucial element is that the background process is irreducible.  Observe that the probability of interest {\it does} depend on the rates $q_{ij}$, and so do the exact asymptotics. We refer to \cite{BTM} for the exact asymptotics of a related infinite-server model; it is noted that the derivation of such precise asymptotics is typically highly involved.}

\RED{\em The above reasoning indicates that the proposed procedure remains valid under more general conditions: the ideas carry over to any situation in which the rates are piecewise constant along the most likely path.
}
\end{remark}

\subsection{Simulation experiments}\label{S44}
We performed experiments featuring a single-node system under Markov modulation. 
In our example the job sizes stem from an exponential distribution. When the background process is in state $i$, the arrival rate is $\lambda_i$, the job-size distribution is exponential with parameter $\mu_i$, and the rate at which the storage level decays is $r_i$, for $i\in\{1,\ldots,d\}.$

The change of measure is then implemented as follows. As pointed out in Section \ref{S42}, per run a path $f$ of the
background process is sampled under the original measure ${\mathbb P}.$ Suppose along this path there are $K$ transitions (remarking that, for compactness, we leave out the argument $f$ here), say at times $t_1$ up to $t_K$; with $t_0=0$ and $t_{K+1}=t$, the state between $t_i$ and $t_{i+1}$ is denoted by $j_i$, for $i=0,\ldots,K.$
Per run a specific change of measure is to be computed, parametrized by the $t_i$ and $j_i$, as follows.

We define
\[P_i(u) := \bar P_i  e^{r_{j_i }u},\:\:\:\:\bar P_i:=   e^{-r_{j_i}t_{i+1}} \prod_{i'=i+1}^K e^{-r_{j_{i'}}(t_{i'+1}-t_{i'})};\]
the product in this expression should be interpreted as $1$ if $i+1>K$. It is readily checked that
\[M(\vt) = \prod_{i=0}^K \exp\left(\lambda_{j_i} \int_{t_i}^{t_{i+1}} \frac{P_i(u)\,\vt}{\mu_{j_i}-P_i(u)\,\vt} {\rm d}u\right).\] Let $\vt\s$ be the maximizing argument of $\vt a-\log M(\vt)$.  

We can now provide the alternative measure ${\mathbb Q}$ for this path of the background process. The number of arrivals between $t_i$ and $t_{i+1}$ (for $i=0,\ldots, K$) becomes Poisson with parameter 
\begin{eqnarray*}\int_{t_i}^{t_{i+1}} \lambda_{j_i} \frac{\mu_{j_i}}{\mu_{j_i}-P_i(u)\,\vt\s}{\rm d}u&=&\frac{\lambda_{j_i}}{r_{j_i}} \log \left( \frac{\mu_{j_i} - \bar P_i e^{r_{j_i} t_i} \vt\s}{\mu_{j_i} e^{-r_{j_i} (t_{i+1}-t_i)} -\bar P_i e^{r_{j_i} t_i} \vt\s} \right)\\
&=&\frac{\lambda_{j_i}}{r_{j_i}} \log \left( \frac{\mu_{j_i} - \bar P_i e^{r_{j_i} t_i} \vt\s}{\mu_{j_i}  -\bar P_i e^{r_{j_i} t_{i+1}} \vt\s} \right)
+\lambda_{j_i}(t_{i+1}-t_i).
\end{eqnarray*}
(where it is noted that this expression is larger than $\lambda_{j_i}(t_{i+1}-t_i)$, which was the parameter under ${\mathbb P}$). 
The density of each of the arrivals between $t_i$ and $t_{i+1}$ becomes
\begin{eqnarray*}\lefteqn{\left(\frac{1}{\mu_{j_i}-P_i(u)\,\vt\s}\right)\left/\int_{t_i}^{t_{i+1}}\left(\frac{1}{\mu_{j_i}-P_i(v)\,\vt\s}\right){\rm d}v\right.}\\
&=&\left(\frac{\mu_{j_i}}{\mu_{j_i}-P_i(u)\,\vt\s}\right)\left/
\frac{1}{r_{j_i}} \log \left( \frac{\mu_{j_i} - \bar P_i e^{r_{j_i} t_i} \vt\s}{\mu_{j_i} e^{-r_{j_i} (t_{i+1}-t_i)} -\bar P_i e^{r_{j_i} t_i} \vt\s} \right)
\right.
\end{eqnarray*}
(rather than a uniform distribution, as was the case under ${\mathbb P}$); sampling from this distribution is easy, since the inverse distribution function can be determined in closed form.
Given an arrival that takes place at time $u$ between $t_i$ and $t_{i+1}$, the job size is exponential with parameter $\mu_{j_i}- P_i(u)\,\vt\s$ (rather than exponential with parameter $\mu_{j_i}$). 

We now describe two examples in which the dimension of the background process is $d=2$,
$q_{12}=q_{21}=2$, and $t=1$.
In the first example we fix $a=3$, ${\bs \lambda}=(2,1)$, ${\bs \mu}=(\frac{1}{2},1),$ and ${\bs r}= (5,1)$, in the second example $a=0.8$, ${\bs \lambda}=(0.9,1)$, ${\bs \mu}=(0.9^{-1},1),$ and ${\bs r}= (0.3,0.6)$. As before, we simulate until the precision of the estimate has reached $\varepsilon=0.1$.
The top two panels in Figs.~\ref{F4}--\ref{F5} should be read as those in Figs.\ \ref{F1}--\ref{F3}; the bottom two panels correspond to the density of the arrival epochs and the rate of the exponential job sizes, respectively, for $f$ the `empirical maximizer' of ${\mathbb I}_f(a)$ (i.e., the maximizer of ${\mathbb I}_f(a)$ over all paths $f$ of the background process that were sampled in the simulation experiment). 

In the first example the thus obtained `optimal path' subsequently visits states 1, 2, and 1, where the
corresponding jump times are $t_1\s=0.654$ and $t_2\s=0.739$, and the decay rate is $0.573$. The mean numbers of arrivals in the three parts of the optimal path are $1.392$, $0.090$ and $0.963$ respectively, whereas for Monte Carlo sampling these are $1.308$, $0.085$ and $0.522$ respectively. 

In the second example the optimal path subsequently visits states 2 and 1, where the
corresponding jump time is $t_1\s=0.790$. In this case the decay rate has the value $0.000806$. The mean numbers of arrivals in the two parts of the optimal path are $0.812$ and $0.195$ respectively, which are slightly higher than the corresponding values under Monte Carlo sampling ($0.790$ and $0.189$ respectively).  Observe that in this example the difference between the two measures is relative small, also reflected by the small value of the decay rate; the event under consideration technically qualifies as `rare' in that $p_n(0.8)\to0$ as $n\to\infty$, but has a relatively high likelihood (e.g.\ as compared to the first example). As a consequence of the fact that both measures almost coincide, the two densities in the bottom-left panel can hardly be distinguished.

We observe that the top panels confirm that in both examples (i) $p_n(a)$ decays roughly exponentially in $n$, (ii) the number of runs needed grows roughly linearly in $n$ (in the first example slightly sublinearly). 

\section{Discussion and concluding remarks}\label{S5}
In this paper we have considered the probability of attaining a value in a rare set $A$ at a fixed point in time $t$: with $A=[a_1,\infty)\times\cdots\times[a_L,\infty)$,
\[p_n(a)={\mathbb P}\left(Y^{(1)}_n(t)\geqslant na_1,\ldots,Y_n^{(L)}(t)\geqslant na_L\right).\] A relevant  related quantity is the probability of having reached the set $A$ {\it before} $t$:
\begin{equation}
\label{PP1}
{\mathbb P}\left(\exists s\leqslant t: Y^{(1)}_n(s)\geqslant na_1,\ldots,Y_n^{(L)}(s)\geqslant na_L\right);\end{equation}
observe that this probability increases to $1$ as $t\to\infty$. Alternatively, one could study the probability that all $a_\ell$ (for $\ell=1,\ldots,L$) are exceeded before $t$, {\it but not necessarily at the same time}:
\begin{equation}
\label{PP2}
{\mathbb P}\left(\exists s_1\leqslant t: Y^{(1)}_n(s_1)\geqslant na_1,\ldots,\exists s_L\leqslant t: Y_n^{(L)}(s_L)\geqslant na_L\right).\end{equation}
Powerful novel sample-path large deviations results by Budhiraja and Nyquist \cite{BN}, which deal with a general class of multi-dimensional shot-noise processes,  may facilitate the development of efficient importance sampling algorithms for non-modulated linear stochastic fluid networks. The results in \cite{BN} do not cover Markov modulation, though.

In the current setup of Section \ref{S4} the speed of the background process is kept fixed, i.e., not scaled by $n$. For modulated diffusions a sample-path large deviation principle has been recently established  in \cite{HU} for the case that the background process is sped up by a factor $n$ (which amounts to multiplying the generator matrix $Q$ by $n$); the rate function decouples into (i)~a part concerning the rare-event behavior of the background process and (ii)~a part concerning the rare-event behavior of the diffusion (conditional on the path of the background process). With a similar result for the Markov-modulated linear stochastic fluid networks that we have studied in this paper, one could potentially set up an efficient importance sampling procedure for the probabilities (\ref{PP1}) and (\ref{PP2}) under this scaling.

\bibliographystyle{plain}

\RED{
\section*{Appendix A}
We here point out how (\ref{SECM}) can be established; the line of reasoning is precisely the same as in the derivation of (\ref{FIRM}) in \cite[Thm. 3.7.4]{DZ}.
First write
\[{\mathbb E}_{\mathbb Q}(L^2I) = {\mathbb E}_{\mathbb Q}( {\rm e}^{-2\vt\s Y_n(t)}{\rm e}^{2n\,\log M(\vt\s)}) = {\rm e}^{-2nI(a)}\,{\mathbb E}_{\mathbb Q}({\rm e}^{-2\vt\s (Y_n(t)-na)} 1_{\{Y_n(t)\geqslant na\}}),\]
which, with $Z_n:=(Y_n(t)-na)/\sqrt{n}$, equals
\[{\rm e}^{-2nI(a)}\,{\mathbb E}_{\mathbb Q}({\rm e}^{-2\vt\s Z_n\sqrt{n}} 1_{\{Z_n\geqslant 0\}}).\]
Observe that ${\mathbb E}_{\mathbb Q}\,Y_n=na$, due to the very choice of ${\mathbb Q}.$ This entails that $Z_n$ converges in distribution to a centered Normal random variable; as can be verified, the corresponding  variance is $\tau$ (where $\tau$ is defined in (\ref{FIRM})). Using the Berry-Esseen-based justification presented in \cite[page 111]{DZ}, we conclude that, as $n\to\infty$,
\[{\mathbb E}_{\mathbb Q}({\rm e}^{-2\vt\s Z_n\sqrt{n}} 1_{\{Z_n\geqslant 0\}})\sim
\int_0^\infty {\rm e}^{-2\vt\s \sqrt{n} \,x} \,\frac{1}{\sqrt{2\pi\tau}} {\rm e}^{-x^2/(2\tau)}{\rm d}x.\]
Completing the square, the right-hand side of the previous display equals, with ${\mathscr N}(\mbox{\sc m},\mbox{\sc v})$ a  normal random variable with mean $\mbox{\sc m}$ and variance $\mbox{\sc v}$,
\[{\rm e}^{2(\vt\s)^2 n\tau} \,{\mathbb P}\big({\mathscr N}(-2\vt\s\sqrt{n}\,\tau,\tau)>0\big)= {\rm e}^{2(\vt\s)^2 n\tau} \,{\mathbb P}\big({\mathscr N}(0,1)>2\vt\s\sqrt{n\tau}\big).\]
Now we use the standard equivalence (as $x\to\infty$)
\[{\mathbb P}({\mathscr N}(0,1)>x) \sim \frac{1}{x}\frac{1}{\sqrt{2\pi}}{\rm e}^{-x^2/2},\]
to obtain
\[
\int_0^\infty {\rm e}^{-2\vt\s \sqrt{n} \,x} \,\frac{1}{\sqrt{2\pi}} {\rm e}^{-x^2/(2\tau)}{\rm d}x\sim
\frac{1}{\sqrt{n}} \,\frac{1}{2\vt\s\sqrt{2\pi\tau}}.\]
Combining the above, we derive the claim:
\[{\mathbb E}_{\mathbb Q}(L^2I)  \sim \frac{1}{\sqrt{n}} \,\frac{1}{2\vt\s\sqrt{2\pi\tau}} {\rm e}^{-2nI(a)}.\]
}

\end{document}